\documentclass[12pt]{article}
\usepackage{amssymb,amsmath,amsthm,tikz,multirow}
\usetikzlibrary{arrows,calc}

\title{Tilings of Sphere by Congruent Pentagons I}
\author{Ka Yue Cheuk, Ho Man Cheung, Min Yan\thanks{Research was supported by Hong Kong RGC General Research Fund 605610 and 606311.} \\ 
Hong Kong University of Science and Technology}

\newtheorem{theorem}{Theorem}
\newtheorem{lemma}[theorem]{Lemma}

\newtheorem{proposition}[theorem]{Proposition}
\newtheorem*{theorem*}{Theorem}

\theoremstyle{definition}
\newtheorem*{definition*}{Definition}
\newtheorem*{case*}{Case}
\newtheorem*{subcase*}{Subcase}

\theoremstyle{remark}

\numberwithin{equation}{section}

\newcommand{\arcThroughThreePoints}[4][]{
\coordinate (middle1) at ($(#2)!.5!(#3)$);
\coordinate (middle2) at ($(#3)!.5!(#4)$);
\coordinate (aux1) at ($(middle1)!1!90:(#3)$);
\coordinate (aux2) at ($(middle2)!1!90:(#4)$);
\coordinate (center) at ($(intersection of middle1--aux1 and middle2--aux2)$);
\draw[#1] 
 let \p1=($(#2)-(center)$),
      \p2=($(#4)-(center)$),
      \n0={veclen(\p1)},       
      \n1={atan2(\y1,\x1)}, 
      \n2={atan2(\y2,\x2)},
      \n3={\n2>\n1?\n2:\n2+360}
    in (#2) arc(\n1:\n3:\n0);
}

\begin{document}

\maketitle

\begin{abstract}
We develop some basic tools for studying edge-to-edge tilings of the sphere by congruent pentagons. Then we prove that there is no edge-to-edge tiling of the sphere by more than $12$ congruent pentagons, under the assumption that the pentagon has edge length combination $a^2b^2c$, $a^3bc$, or $a^3b^2$ ($a,b,c$ distinct), and there is a tile with all vertices having degree $3$. 
\end{abstract}

\section{Introduction}

Mathematicians have studied tilings for more than 100 years. A lot is known about tilings of the plane or the Euclidean space. However, results about tilings of the sphere are relatively rare. A major achievement in this regard is the complete classification of edge-to-edge tilings of the sphere by congruent triangles \cite{so,ua}. For tilings of the sphere by congruent pentagons, we completely classified the minimal case of $12$ tiles \cite{ay1,gsy}.

The spherical tilings should be easier to study than the planar tilings, simply because the former involves only finitely many tiles. The classifications in \cite{gsy,ua} not only give the complete list of tiles, but also the ways the tiles are fit together. It is not surprising that such kind of classifications for the planer tilings are only possible under various symmetry conditions, because the quotients of the plane by the symmetries often become compact.

Like the earlier works, we restrict ourselves to edge-to-edge tilings of the sphere by congruent polygons, such that all vertices have degree $\ge 3$. These are mild and natural assumptions that simplify the discussion. The polygon in such a tiling must be triangle, quadrilateral, or pentagon \cite{ua2}. We believe that pentagonal tilings should be relatively easier to study than quadrilateral ones because $5$ is an ``extreme'' among $3$, $4$, $5$. Indeed, in Section \ref{vet}, we find various restrictions on pentagonal tilings of the sphere. Similar restrictions for quadrilateral tilings are much weaker, which is the primary reason that the study of quadrilateral tilings is much harder \cite{ua2}. Our classification program starts with Proposition \ref{base_tile}, which says that a pentagonal tilings of the sphere must have a tile, such that four vertices have degree $3$, and the fifth vertex has degree $3$, $4$ or $5$. We call such a special tile $3^5$-, $3^44$-, or $3^45$-tile. Our strategy is to first try to tile the neighborhood of this special tile. The neighborhood tiling gives us a lot of information on edges and angles of the pentagon. Then we use the information to determine the {\em anglewise vertex combination} (all the possible angle combinations at vertices, abbreviated as AVC). The AVC guides us to further construct the tiling beyond the neighborhood, and we eventually construct the tiling of the sphere.

Proposition \ref{base_tile} and the subsequent Propositions \ref{base_tile2} and \ref{base_tile3} give only the combinatorial structure of the special tile. Further information on edge lengths is provided by \cite[Proposition 7]{gsy} (and its extension Proposition \ref{edge_combo} in this paper), which says that the pentagon can only have $5$ possible edge length combinations ($a,b,c$ are distinct): $a^5$ (equilateral), $a^4b$ (almost equilateral), $a^3b^2$, $a^3bc$, $a^2b^2c$ (variable edge length). The main result of this paper is the classification of the case of variable edge length, under the additional assumption that there is a $3^5$-tile.

\begin{theorem*}
There is no edge-to-edge spherical tiling by more than $12$ congruent pentagons, such that there is a tile with all vertices having degree $3$, and the pentagon has edge length combination $a^2b^2c$, $a^3bc$, or $a^3b^2$ ($a,b,c$ distinct).
\end{theorem*} 

This is the first of the following series of papers on the classification of edge-to-edge tilings the sphere by congruent pentagons.
\begin{enumerate}
  \item[O.] The minimal case of tilings by $12$ congruent pentagons. \cite{gsy}
  \item[I.] Edge length combinations $a^3b^2$, $a^3bc$, $a^2b^2c$, and there is a $3^5$-tile. [this paper]
  \item[II.] Edge length combinations $a^3bc$, $a^2b^2c$, and there is no $3^5$-tile. \cite{wy1}
  \item[III.] Edge length combination $a^3b^2$, and there is no $3^5$-tile.  \cite{wy2}
  \item[IV.] Equilateral length combination $a^5$.  \cite{ay2}
\end{enumerate}
After the ``0-th'' paper \cite{gsy}, we may always assume more than $12$ tiles. After the five papers, the remaining case is the almost equilateral  length combination $a^4b$. The case is the most complicated, and allows many tilings not appearing in the other cases.

Sections \ref{vet}, \ref{constraint1}, and \ref{constraint2} are general results. Section \ref{vet} is the basic results on the distribution of vertices, angles, and edges. Section \ref{constraint1} gives a crude yet very effective geometrical constraint on pentagon with two pairs of equal edges. Section \ref{constraint2} gives precise constraints on spherical pentagons in terms of spherical trigonometry. The constraints are used to eliminate two cases in this paper. 

Sections \ref{nd} and \ref{3a2b_tiling} classify the specific tilings of this paper, which assumes variable edge length and existence of a $3^5$-tile. Section \ref{nd} studies the neighborhood tiling of a $3^5$-tile. In Theorems \ref{thm1} and \ref{thm2}, the information from the neighborhood tiling quickly implies that the edge length combinations $a^2b^2c$ and $a^3bc$ do not admit tilings by more than $12$ congruent pentagons. For the edge length combination $a^3b^2$, Proposition \ref{main3} shows that there are $4$ possible neighborhood tilings. In Section \ref{3a2b_tiling}, we further prove that none of the $4$ possible neighborhood tilings lead to tilings of the sphere.

\section{Vertex, Angle and Edge}
\label{vet}

\subsubsection*{Vertex}

Consider an edge-to-edge tiling of the sphere by pentagons, such that all vertices have degree $\ge 3$. Let $v,e,f$ be the numbers of vertices, edges, and tiles. Let $v_k$ be the number of vertices of degree $k$. We have
\begin{align*}
2
&=v-e+f, \\
2e
&=5f=\sum_{k=3}^{\infty}kv_k=3v_3+4v_4+5v_5+\cdots, \\
v
&=\sum_{k=3}^{\infty}v_k=v_3+v_4+v_5+\cdots.
\end{align*}
Then it is easy to derive $2v=3f+4$ and  
\begin{align}
\frac{f}{2}-6
&=\sum_{k\ge 4}(k-3)v_k=v_4+2v_5+3v_6+\cdots, \label{vcountf} \\
v_3
&=\sum_{k\ge 4}(3k-10)v_k=2v_4+5v_5+8v_6+\cdots. \label{vcountv}
\end{align}
By \eqref{vcountf}, $f$ is an even integer $\ge 12$. Since tilings by $12$ congruent pentagons have been classified by \cite{ay1, gsy}, we may assume $f>12$. We also note that by \eqref{vcountf}, $f=14$ implies $v_4=1$ and $v_5=v_6=\cdots=0$. By \cite[Theorem 1]{yan}, this is impossible. Therefore we will always assume $f$ is an even integer $\ge 16$.

By \eqref{vcountv}, most vertices have degree $3$. We call vertices of degree $>3$ {\em high degree} vertices.

\begin{proposition}\label{base_tile}
Any pentagonal tiling of the sphere has a tile, such that four vertices have degree $3$ and the fifth vertex has degree $3$, $4$ or $5$.
\end{proposition}

We call three types of tiles in the proposition $3^5$-tile, $3^44$-tile, and $3^45$-tile. 

\begin{proof}
Suppose the tile described in the proposition does not exist. Then any tile either has at least one vertex of degree $\ge 6$, or has at least two vertices of degree $4$ or $5$. Since a degree $k$ vertex is shared by at most $k$ tiles, the number of tiles of first kind is $\le \sum_{k\ge 6}kv_k$, and the number of tiles of the second kind is $\le\frac{1}{2}(4v_4+5v_5)$. Therefore we have
\[
f\le 2v_4+\frac{5}{2}v_5+\sum_{k\ge 6}kv_k.
\]
On the other hand, by \eqref{vcountf}, we have
\[
f-\left(2v_4+\frac{5}{2}v_5+\sum_{k\ge 6}kv_k \right)
= 12 + \frac{1}{2}v_5+\sum_{k\ge 4}(k-6)v_k
>0.
\]
We get a contradiction.
\end{proof}

\begin{proposition}\label{base_tile2}
If a pentagonal tiling of the sphere has no $3^5$-tile, then $f\ge 24$. Moreover, if $f=24$, then each tile is a $3^44$-tile. 
\end{proposition}

\begin{proof}
Suppose there is no $3^5$-tile. Then any tile has at least one vertex of degree $\ge 4$. Therefore we have $f\le\sum_{k\ge 4}kv_k$. By \eqref{vcountf}, we further have
\begin{align*}
f=2f-f
&\ge 24+ \sum_{k\ge 4}4(k-3)v_k -\sum_{k\ge 4}kv_k \\
&=24+ \sum_{k\ge 4}3(k-4)v_k
\ge 24.
\end{align*}
Moreover, the equality happens if and only if $v_5=v_6=\cdots=0$ and $f= 4v_4$. Since there is no $3^5$-tile, this means that each tile is a $3^44$-tile.
\end{proof}

\begin{proposition}\label{base_tile3}
If a pentagonal tiling of the sphere has no $3^5$-tile and $3^44$-tile, then $f\ge 60$. Moreover, if $f=60$, then each tile is a $3^45$-tile. 
\end{proposition}

\begin{proof}
Suppose there is no $3^5$-tile and $3^44$-tile. Then any tile either has at least one vertex of degree $\ge 5$, or has at least two vertices of degree $4$. Therefore we have $f\le \frac{1}{2}4v_4+\sum_{k\ge 5}kv_k$. By \eqref{vcountf}, we further have
\begin{align*}
f=5f-4f
&\ge 60+\sum_{k\ge 4}5(k-3)v_k-8v_4-\sum_{k\ge 5}4kv_k \\
&=60+2v_4+\sum_{k\ge 5}5(k-6)v_k
\ge 60.
\end{align*}
Moreover, the equality happens if and only if $v_4=v_6=\cdots=0$ and $f= 5v_4$. Since there is no $3^5$-tile and $3^44$-tile, this means that each tile is a $3^45$-tile.
\end{proof}

\subsubsection*{Angle}

The most basic property about angles is that the sum of all angles ({\em angle sum}) at a vertex is $2\pi$. Another basic property is the sum of all angles in the pentagon.

\begin{proposition}\label{anglesum}
If all tiles in a tiling of sphere by $f$ pentagons have the same five angles $\alpha,\beta,\gamma,\delta,\epsilon$, then 
\begin{equation}\label{asump}
\alpha+\beta+\gamma+\delta+\epsilon
=3\pi + \frac{4}{f}\pi.
\end{equation}
\end{proposition}

\begin{proof}
Since the angle sum at each vertex is $2\pi$, the total sum of all angles is $2\pi v$. Since all tiles are (angle) congruent, the sum $\Sigma$ of five angles is the same for all the tiles, and the total sum of all angles is also $f\Sigma$. Therefore we have $2\pi v = f\Sigma$. Then by $3f=2v-4$, we get
\[
\Sigma = 2\pi\frac{v}{f} = 3\pi + \frac{4}{f}\pi. \qedhere
\]
\end{proof}

The proposition does not require that the angles are arranged in the same way in all tiles. Moreover, if we additionally know that all edges are straight (i.e., great arcs), then all tiles have the same area $\Sigma-3\pi$, and \eqref{asump} follows from the fact that the total area $(\Sigma-3\pi)f$ is the area $4\pi$ of the sphere. The proposition does not require that the edges are straight.

The angles in Proposition \ref{anglesum} refer to the values, and some angles among the five may have the same value. For example, if the five values are $\alpha,\alpha,\alpha,\beta,\beta$, with $\alpha\ne\beta$ (different values), then we say the pentagon has {\em angle combination} $\alpha^3\beta^2$. The following is about the distribution of angle values.

\begin{proposition}\label{deg3a}
If an angle appears at every degree $3$ vertex in a tiling of sphere by pentagons with the same angle combinations, then the angle appears at least $2$ times in the pentagon.
\end{proposition}

\begin{proof}
If an angle $\theta$ appears only once in the pentagon, then the total number of times $\theta$ appears in the whole tiling is $f$, and the total number of non-$\theta$ angles is $4f$. If we also know that $\theta$ appears at every degree $3$ vertex, then $f\ge v_3$, and non-$\theta$ angles appear $\le 2v_3$ times at degree $3$ vertices. Moreover, non-$\theta$ angles appear $\le \sum_{k\ge 4}kv_k$ times at high degree vertices. Therefore 
\[
4v_3 \le 4f \le 2v_3+\sum_{k\ge 4}kv_k.
\]
On the other hand, by \eqref{vcountv}, we have
\[
4v_3-\left( 2v_3+\sum_{k\ge 4}kv_k \right)
=\sum_{k\ge 4}[2(3k-10)-k]v_k
=\sum_{k\ge 4}5(k-4)v_k
\ge 0.
\]
Then we get $v_5=v_6=\cdots=0$ and $f=v_3=2v_4$. This contradicts \eqref{vcountf}.
\end{proof}

Unlike Proposition \ref{anglesum}, which is explicitly about the values of angles, Proposition \ref{deg3a} only counts the number of angles. The key in counting is to distinguish angles. We may use the value as the criterion for the two angles to be the ``same''. We may also use the edge lengths bounding the angles as the criterion. The observation will be used in the proof of Proposition \ref{edge_combo}. The observation also applies to the subsequent Propositions \ref{deg3b}, \ref{deg3c}, \ref{hdeg}.

\begin{proposition}\label{deg3b}
If an angle appears at least twice at every degree $3$ vertex in a tiling of sphere by pentagons with the same angle combination, then the angle appears at least $3$ times in the pentagon. 
\end{proposition}

\begin{proof}
If an angle $\theta$ appears only once in the pentagon, then by Proposition \ref{deg3a}, it cannot appear at every degree $3$ vertex.

If an angle $\theta$ appears twice in the pentagon, then the total number of times $\theta$ appears in the whole tiling is $2f$, and the total number of non-$\theta$ angles is $3f$. If we also know that $\theta$ appears at least twice at every degree $3$ vertex, then $2f\ge 2v_3$, and non-$\theta$ angles appear $\le v_3$ times at degree $3$ vertices. Moreover, non-$\theta$ angles appear $\le \sum_{k\ge 4}kv_k$ times at high degree vertices. Therefore 
\[
3v_3 \le 3f \le v_3+\sum_{k\ge 4}kv_k.
\]
This leads to the same contradiction as in the proof of Proposition \ref{deg3a}.
\end{proof}

The proof of Proposition \ref{deg3b} can be easily modified to get the following.

\begin{proposition}\label{deg3c}
If two angles together appear at least twice at every degree $3$ vertex in a tiling of sphere by pentagons with the same angle combination, then the two angles together appear at least $3$ times in the pentagon. 
\end{proposition}

The following is about angles not appearing at degree $3$ vertices.

\begin{proposition}\label{hdeg}
Suppose an angle $\theta$ does not appear at degree $3$ vertices in a tiling of sphere by pentagons with the same angle combination.
\begin{enumerate}
\item There can be at most one such angle $\theta$.
\item The angle $\theta$ appears only once in the pentagon.
\item $2v_4+v_5\ge 12$.
\item One of $\alpha\theta^3$, $\theta^4$, $\theta^5$ is a vertex, where $\alpha\ne\theta$.
\end{enumerate}
\end{proposition}

The first statement implies that the angle $\alpha$ in the fourth statement must appear at a degree $3$ vertex.

\begin{proof}
Suppose two angles $\theta_1$ and $\theta_2$ do not appear at degree $3$ vertices. Then the total number of times these two angles appear is at least $2f$, and is at most the total number $\sum_{k\ge 4}kv_k$ of angles at high degree vertices. Therefore we have $2f\le \sum_{k\ge 4}kv_k$. On the other hand, by \eqref{vcountf}, we have
\[
2f - \sum_{k\ge 4}kv_k  
= 24 + \sum_{k\ge 4}3(k-4)v_k
>0.
\]
The contradiction proves the first statement.

The argument above also applies to the case $\theta_1=\theta_2$, which means the same angle appearing at least twice in the pentagon. This proves  the second statement. 

The first two statements imply that $\theta$ appears exactly $f$ times. Since this should be no more than the total number $\sum_{k\ge 4}kv_k$ of angles at high degree vertices, by \eqref{vcountf}, we have
\[
0
\ge f - \sum_{k\ge 4}kv_k  
= 12 - 2v_4 - v_5 + \sum_{k\ge 6}(k-6)v_k.
\]
This implies the third statement.

For the last statement, we assume that $\alpha\theta^3,\theta^4,\theta^5$ are not vertices. This means that $\theta$ appears at most twice at any degree $4$ vertex, and at most four times at any degree $5$ vertex. Since $\theta$ also does not appear at degree $3$ vertices, the total number of times $\theta$ appears is $\le 2v_4+4v_5+\sum_{k\ge 6}kv_k$. However, the number of times $\theta$ appears should also be $f$. Therefore $f\le 2v_4+4v_5+\sum_{k\ge 6}kv_k$. On the other hand, by \eqref{vcountf}, we have
\[
f-\left(2v_4+4v_5+\sum_{k\ge 6}kv_k \right)
= 12 + \sum_{k\ge 6}(k-6)v_k
>0.
\]
We get a contradiction.
\end{proof}

\subsubsection*{Edge}

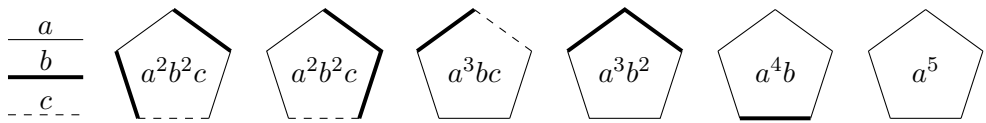
\begin{figure}[htp]
\centering
\begin{tikzpicture}[>=latex,scale=1]

\begin{scope}[xshift=-2.2cm]

\draw
	(0,0.4) -- node[above=-2] {\small $a$} ++(1,0);

\draw[line width=1.5]
	(0,-0.1) -- node[above=-2] {\small $b$} ++(1,0);

\draw[dashed]
	(0,-0.6) -- node[above=-2] {\small $c$} ++(1,0);

\end{scope}


\draw
	(90:0.8) -- (162:0.8)
	(18:0.8) -- (-54:0.8);

\draw[line width=1.5]
	(90:0.8) -- (18:0.8)
	(162:0.8) -- (234:0.8);

\draw[dashed]
	(-54:0.8) -- (-126:0.8);
	
\node at (0,0) {\small $a^2b^2c$};


\begin{scope}[xshift=2cm]

\draw
	(90:0.8) -- (162:0.8) -- (234:0.8);

\draw[line width=1.5]
	(90:0.8) -- (18:0.8) -- (-54:0.8);

\draw[dashed]
	(-54:0.8) -- (-126:0.8);
	
\node at (0,0) {\small $a^2b^2c$};

\end{scope}


\begin{scope}[xshift=4cm]

\foreach \x in {2,...,4}
\draw[rotate=72*\x]
	(90:0.8) -- (18:0.8);

\draw[line width=1.5]
	(162:0.8) -- (90:0.8);

\draw[dashed]
	(18:0.8) -- (90:0.8);
	
\node at (0,0) {\small $a^3bc$};

\end{scope}


\begin{scope}[xshift=6cm]

\draw
	(162:0.8) -- (234:0.8) -- (-54:0.8) -- (18:0.8);

\draw[line width=1.5]
	(162:0.8) -- (90:0.8) -- (18:0.8);
	
\node at (0,0) {\small $a^3b^2$};

\end{scope}


\begin{scope}[xshift=8cm]

\foreach \x in {-1,...,2}
\draw[rotate=72*\x]
	(90:0.8) -- (18:0.8);

\draw[line width=1.5]
	(-54:0.8) -- (-126:0.8);
	
\node at (0,0) {\small $a^4b$};

\end{scope}


\begin{scope}[xshift=10cm]

\foreach \x in {1,...,5}
\draw[rotate=72*\x]
	(90:0.8) -- (18:0.8);
\node at (0,0) {\small $a^5$};

\end{scope}

\end{tikzpicture}
\caption{Edges in the pentagon suitable for tiling, $a,b,c,d,e$ distinct.}
\label{edges1}
\end{figure}

Two pentagons have the same edge length {\em combination} if the five edge lengths are equal. For example, the first two pentagons in Figure \ref{edges1} have the same edge length combination $a^2b^2c$. If we further have the same edge length {\em arrangement} in the two pentagons, then the two pentagons are {\em edge congruent}. For example, the first two pentagons in Figure \ref{edges1} have the respective edge length arrangements $a,b,a,b,c$ and $a,a,b,b,c$, and therefore they are not edge congruent.

The subsequent discussion does not have any requirement on the values of angles.

\begin{proposition}\label{edge_combo}
In an edge-to-edge tiling of the sphere by edge congruent pentagons, the edge lengths of the pentagon is arranged in one of the six ways in Figure \ref{edges1}. 
\end{proposition}

\begin{proof}
The proposition is an extension of \cite[Proposition 7]{gsy}. We first follow the argument in the earlier paper.

By purely numerical consideration, there are seven possible edge length combinations
\[
abcde,\;
a^2bcd,\;
a^2b^2c,\;
a^3bc,\;
a^3b^2,\;
a^4b,\;
a^5.
\]
For $abcde$, without loss of generality, we may assume that the edges are arranged as in the first of Figure \ref{edges2}. By Proposition \ref{base_tile}, one of the two bottom vertices (shared by $b,c$, and shared by $c,d$) has degree $3$. Without loss of generality, we may assume that the vertex shared by $c,d$ has degree $3$. Let $x$ be the third edge at the vertex. Then $x,c$ are adjacent in a tile, and $x,d$ are adjacent in another tile. Since there is no edge in the pentagon that is adjacent to both $c$ and $d$, we get a contradiction. 

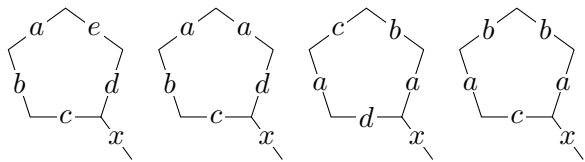
\begin{figure}[htp]
\centering
\begin{tikzpicture}[>=latex,scale=1]


\draw
	(18:0.8) -- node[fill=white,inner sep=1] {\small $e$}
	(90:0.8) -- node[fill=white,inner sep=1] {\small $a$} 
	(162:0.8) -- node[fill=white,inner sep=1] {\small $b$} 
	(234:0.8) -- node[fill=white,inner sep=1] {\small $c$} 
	(-54:0.8) -- node[fill=white,inner sep=1] {\small $d$} 
	(18:0.8)
	(-54:0.8) -- node[fill=white,inner sep=1] {\small $x$} 
	(-54:1.5);


\begin{scope}[xshift=2cm]

\draw
	(18:0.8) -- node[fill=white,inner sep=1] {\small $a$}
	(90:0.8) -- node[fill=white,inner sep=1] {\small $a$} 
	(162:0.8) -- node[fill=white,inner sep=1] {\small $b$} 
	(234:0.8) -- node[fill=white,inner sep=1] {\small $c$} 
	(-54:0.8) -- node[fill=white,inner sep=1] {\small $d$} 
	(18:0.8)
	(-54:0.8) -- node[fill=white,inner sep=1] {\small $x$} 
	(-54:1.5);

\end{scope}


\begin{scope}[xshift=4cm]

\draw
	(18:0.8) -- node[fill=white,inner sep=1] {\small $b$}
	(90:0.8) -- node[fill=white,inner sep=1] {\small $c$} 
	(162:0.8) -- node[fill=white,inner sep=1] {\small $a$} 
	(234:0.8) -- node[fill=white,inner sep=1] {\small $d$} 
	(-54:0.8) -- node[fill=white,inner sep=1] {\small $a$} 
	(18:0.8)
	(-54:0.8) -- node[fill=white,inner sep=1] {\small $x$} 
	(-54:1.5);
	
\end{scope}


\begin{scope}[xshift=6cm]

\draw
	(18:0.8) -- node[fill=white,inner sep=1] {\small $b$}
	(90:0.8) -- node[fill=white,inner sep=1] {\small $b$} 
	(162:0.8) -- node[fill=white,inner sep=1] {\small $a$} 
	(234:0.8) -- node[fill=white,inner sep=1] {\small $c$} 
	(-54:0.8) -- node[fill=white,inner sep=1] {\small $a$} 
	(18:0.8)
	(-54:0.8) -- node[fill=white,inner sep=1] {\small $x$} 
	(-54:1.5);

\end{scope}

\end{tikzpicture}
\caption{Not suitable for tiling.}
\label{edges2}
\end{figure}

The combination $a^2bcd$ has two possible arrangements, illustrated (without loss of generality) in the second (adjacent $a$) and third (separated $a$) of Figure \ref{edges2}. Similar to the case of $abcde$, we may assume a degree $3$ vertex with the third edge $x$. In the second of Figure \ref{edges2}, the edge $x$ is adjacent to both $c$ and $d$, a contradiction. In the third of Figure \ref{edges2}, the edge $x$ is adjacent to both $a$ and $d$, a contradiction. 

The combination $a^2b^2c$ has three possible arrangements, illustrated in the first and second of Figure \ref{edges1} and fourth of Figure \ref{edges2}. In the fourth of Figure \ref{edges2}, we may assume a degree $3$ vertex with the third edge $x$. The edge $x$ is adjacent to both $a$ and $c$, a contradiction.

The combination $a^3bc$ has two possible arrangements, illustrated in the third of Figure \ref{edges1} ($b,c$ adjacent) and left of Figure \ref{edges3} ($b,c$ separated). In case $b,c$ are separated, the pentagon has one $a^2$-angle denote by $\alpha$, two $ab$-angles, and two $ac$-angles. Since there are no $b^2$-angle, $c^2$-angle, $bc$-angle, any degree $3$ vertex must be one of three on the left of Figure \ref{edges3}. This implies that the $a^2$-angle $\alpha$ appears at every degree $3$ vertex. By Proposition \ref{deg3a} (and the remark after the proof of the proposition), the pentagon should have at least two $a^2$-angles, a contradiction.

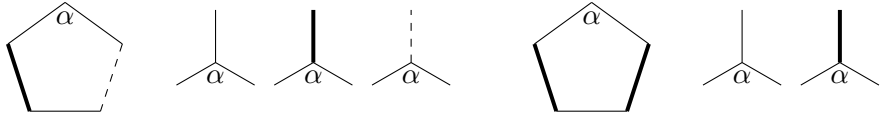
\begin{figure}[htp]
\centering
\begin{tikzpicture}[>=latex,scale=1]


\foreach \x in {0,1,-2}
\draw[rotate=72*\x]
	(90:0.8) -- (18:0.8);

\draw[line width=1.5]
	(162:0.8) -- (234:0.8);

\draw[dashed]
	(18:0.8) -- (-54:0.8);
	
\node at (90:0.6) {\small $\alpha$};

\foreach \a in {0,1,2}
\draw[xshift=2 cm+1.3*\a cm]
	(-30:0.6) -- (0,0) node[below=-1] {\small $\alpha$} -- (210:0.6);
	
\draw
	(2,0) -- ++(0,0.7);
	
\draw[line width=1.5]
	(3.3,0) -- ++(0,0.7);

\draw[dashed]
	(4.6,0) -- ++(0,0.7);
	

\begin{scope}[xshift=7cm]

\foreach \x in {0,1,-2}
\draw[rotate=72*\x]
	(90:0.8) -- (18:0.8);

\draw[line width=1.5]
	(162:0.8) -- (234:0.8)
	(18:0.8) -- (-54:0.8);

\node at (90:0.6) {\small $\alpha$};

\foreach \a in {0,1}
\draw[xshift=2 cm+1.3*\a cm]
	(-30:0.6) -- (0,0) node[below=-1] {\small $\alpha$} -- (210:0.6);
	
\draw
	(2,0) -- ++(0,0.7);
	
\draw[line width=1.5]
	(3.3,0) -- ++(0,0.7);
	
\end{scope}

\end{tikzpicture}
\caption{Also not suitable for tiling.}
\label{edges3}
\end{figure}

The combination $a^3b^2$ has two possible arrangements, illustrated in the fourth of Figure \ref{edges1} ($b$ adjacent) and right of Figure \ref{edges3} ($b$ separated). In case the two $b$-edges are separated, the pentagon has one $a^2$-angle $\alpha$ and four $ab$-angles. Since there is no $b^2$-angle, any degree $3$ vertex must be one of two on the right of Figure \ref{edges3}. This implies that the $a^2$-angle $\alpha$ appears at every degree $3$ vertex. By Proposition \ref{deg3a}, the pentagon should have at least two $a^2$-angles, a contradiction.
\end{proof}

The following gives further details about the special tile in Proposition \ref{base_tile} in case of the first edge length arrangement in Figure \ref{edges1}.

\begin{proposition}\label{edge_combo2}
An edge-to-edge tiling of the sphere by pentagons of edge length arrangement in the first of Figure \ref{edges1} cannot be a $3^5$-tile. Moreover, if the pentagon is a $3^44$- or $3^45$-tile, then the fifth vertex of degree $4$ or $5$ is opposite to the $c$-edge. 
\end{proposition}

\begin{proof}
The edge length arrangement in the first of Figure \ref{edges1} has three $ab$-vertices (shared by $a$-edge and $b$-edge). If any such $ab$-vertex has degree $3$, then the third edge at the vertex is adjacent to both $a$ and $b$. The only such edge in the pentagon is $c$. 

If two $ab$-vertices of degree $3$ are adjacent, then we have two $c$-edges at the two vertices, and these two $c$-edges belong to the same pentagon, a contradiction. Therefore we cannot have adjacent $ab$-vertices of degree $3$. The observation implies the conclusion of the proposition.
\end{proof}

\section{Non-symmetric Pentagon}
\label{constraint1}

\begin{figure}[htp]
\centering
\begin{tikzpicture}

\draw
	(0,0.7) -- node[fill=white,inner sep=1] {\small $b$}
	(-1,0) -- node[fill=white,inner sep=1] {\small $a$}
	(-0.7,-1) -- 
	(0.7,-1) -- node[fill=white,inner sep=1] {\small $a$}
	(1,0) -- node[fill=white,inner sep=1] {\small $b$}
	cycle;

\node at (0,0.5) {\small $\alpha$};	
\node at (-0.75,-0.05) {\small $\beta$};
\node at (0.8,-0.1) {\small $\gamma$};	
\node at (-0.6,-0.8) {\small $\delta$};	
\node at (0.6,-0.8) {\small $\epsilon$};	

\end{tikzpicture}
\caption{Geometrical constraint for pentagon.}
\label{geom1}
\end{figure}
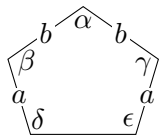

Consider the pentagon in Figure \ref{geom1}, with two $b$-edges and two $a$-edges as indicated. By \cite[Lemma 21]{gsy}, we have $\beta=\gamma$ if and only if $\delta=\epsilon$. Of course the equalities mean that the pentagon is symmetric. We try to prove the following stronger version of the result.

\begin{lemma}\label{geometry1}
If the spherical pentagon in Figure \ref{geom1} is simple and has two pairs of equal edges $a$ and $b$, then $\beta>\gamma$ is equivalent to $\delta<\epsilon$.
\end{lemma}

By simple polygon, we mean the boundary does not intersect itself. By \cite[Lemma 1]{gsy}, in an edge-to-edge tiling of the sphere by congruent pentagons, the pentagon must be simple.

To simplify the discussion, we will always use strict inequalities in this section. For example, this means that the opposite of $\alpha>\pi$ will be $\alpha<\pi$. The special case of equalities can be easily analysed.

The concept of the interior of a simple polygon is relative on sphere. Once we designate ``inside'', then we have the ``outside''. Moreover, the inside is the ``outside of the outside''. We note that the reformulation of Lemma \ref{geometry1} in terms of the outside is equivalent to the original formulation in terms of the inside. We will also use the sine law and the following well known result in the spherical trigonometry: If a spherical triangle has the angles $\alpha,\beta,\gamma<\pi$ and has the edges $a,b,c$ opposite to the angles, then $\alpha>\beta$ if and only if $a>b$.

\begin{figure}[htp]
\centering
\begin{tikzpicture}


\draw
	(0,1) -- (-1.2,-1) -- (-0.3,-1.5) -- (0.3,-0.5) -- (1.2,-1) -- cycle;
\draw[dashed]
	(1.2,-1) -- (-1.2,-1)
	(0.3,-0.5) -- ++(-0.9,0.5);

\node at (0,1.2) {\small $A$};
\node at (-1.4,-1) {\small $B$};
\node at (1.4,-1) {\small $C$};
\node at (-0.3,-1.7) {\small $D$};
\node at (0.05,-0.55) {\small $E$};
\node at (0.1,-1.2) {\small $F$};
\node at (-0.8,0) {\small $G$};

\node[fill=white,inner sep=1] at (0.4,0.3) {\small $b$};
\node[fill=white,inner sep=1] at (-0.4,0.3) {\small $b$};
\node[fill=white,inner sep=1] at (-0.7,-1.3) {\small $a$};
\node[fill=white,inner sep=1] at (0.5,-0.6) {\small $a$};

\node at (0,0.7) {\small $\alpha$};
\node at (-0.9,-0.9) {\small $\beta$};
\node at (0.75,-0.6) {\small $\gamma$};
\node at (-0.35,-1.25) {\small $\delta$};
\node at (0.3,-0.35) {\small $\epsilon$};


\begin{scope}[xshift=3.5cm]

\draw
	(0,1) -- (-1.2,-0.2) -- (-0.8,-1.5) -- (0.8,-1.5) -- (1.2,-0.2) -- cycle;
\draw[dashed]
	(1.2,-0.2) -- (-1.2,-0.2);

\node at (0,1.2) {\small $A$};
\node at (-1.4,-0.2) {\small $B$};
\node at (1.4,-0.2) {\small $C$};
\node at (-0.9,-1.65) {\small $D$};
\node at (0.9,-1.65) {\small $E$};

\node[fill=white,inner sep=1] at (-1,-0.9) {\small $a$};
\node[fill=white,inner sep=1] at (1,-0.9) {\small $a$};
\node[fill=white,inner sep=1] at (-0.5,0.5) {\small $b$};
\node[fill=white,inner sep=1] at (0.5,0.5) {\small $b$};

\node at (0,0.7) {\small $\alpha$};  
\node at (-0.65,-1.3) {\small $\delta$};
\node at (0.7,-1.3) {\small $\epsilon$};  
\node at (-0.9,-0.4) {\small $\beta'$};
\node at (0.95,-0.4) {\small $\gamma'$};

\end{scope}


\begin{scope}[xshift=7cm]

\draw
	(0,1) -- (-1.2,-1.5) -- (-0.35,-0.7) -- (0.35,-0.7) -- (1.2,-1.5) -- cycle;
\draw[dashed]
	(1.2,-1.5) -- (-1.2,-1.5);
	
\node at (0,1.2) {\small $A$};
\node at (-1.4,-1.5) {\small $B$};
\node at (1.4,-1.5) {\small $C$};

\node at (0,0.7) {\small $\alpha$};
\node at (-0.8,-0.95) {\small $\beta$};
\node at (0.75,-0.9) {\small $\gamma$};
\node at (-0.3,-0.5) {\small $\delta$};
\node at (0.3,-0.55) {\small $\epsilon$};
		
\node at (-0.7,-1.3) {\small $\beta'$};
\node at (0.6,-1.25) {\small $\gamma'$};
\node at (-0.3,-0.9) {\small $\delta'$};
\node at (0.25,-0.9) {\small $\epsilon'$};

\node[fill=white,inner sep=1] at (0.5,0) {\small $b$};
\node[fill=white,inner sep=1] at (-0.5,0) {\small $b$};

\end{scope}

\end{tikzpicture}
\caption{Constraint for pentagon and for quadrilateral.}
\label{geom2}
\end{figure}
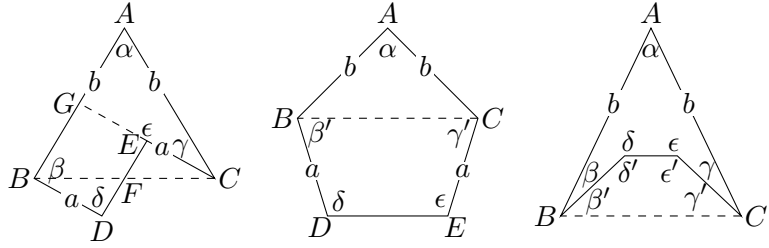

Let $A,B,C,D,E$ be the vertices at the angles $\alpha,\beta,\gamma,\delta,\epsilon$. Among the two great arcs connecting $B$ and $C$, take the one with length $<\pi$ to form the edge $BC$. The edge $BC$ is indicated by the dashed line in Figure \ref{geom2}. Since $AB$ and $AC$ intersect only at one point $A$ and have the same length $b$, we get $b<\pi$. Since all three edges $AB,AC,BC$ have lengths $<\pi$, by the sine law, among the two triangles bounded by the three edges, one has all three angles $<\pi$. We denote this triangle by $\triangle ABC$. The angle $\angle BAC$ of $\triangle ABC$ is either $\alpha$ of the pentagon, or its complement $2\pi-\alpha$. Since the inside and outside versions of the lemma are equivalent, we will always assume $\alpha=\angle BAC<\pi$ in the subsequent discussion.

The pentagon is obtained by choosing $D,E$, and then connecting $B$ to $D$, $C$ to $E$, and $D$ to $E$ by great arcs. Since $BC<\pi$, we find that $BC$ does not intersect $BD$ and $CE$, and the intersection of $BC$ and $DE$ is at most one point. 

If $BC$ and $DE$ intersect at one point $F$, then one of $D,E$ is inside $\triangle ABC$ and one is outside (we omit the ``equality case'' of $D$ or $E$ is on $BC$). The first of Figure \ref{geom2} shows the case $D$ is outside and $E$ is inside. Since $BC<\pi$, the interiors of $BD$ and $CE$ do not intersect $BC$. This implies that 
\[
\beta>\angle ABC=\angle ACB> \gamma.
\]
On the other hand, since $AC=b<\pi$ and $BC<\pi$, the prolongation of $CE$ intersects the boundary of $\triangle ABC$ at a point $G$ on $AB$. Using $AB<\pi$, $\alpha<\pi$, $\gamma<\angle ACB<\pi$ and applying the sine law to $\triangle ACG$, we get $CG<\pi$, so that $a=CE<CG<\pi$. Using $a<\pi$, $BF<BC<\pi$, $CF<BC<\pi$, $\angle BFD=\angle CFE<\pi$ and applying the sine law to $\triangle BDF$ and $\triangle CEF$, we get $\angle BDF<\pi$ and $\angle CEF<\pi$. Therefore
\[
\delta=\angle BDF<\pi<2\pi-\angle CEF=\epsilon.
\] 
This proves that $D$ outside and $E$ inside imply $\beta>\gamma$ and $\delta<\epsilon$. Similarly, $D$ inside and $E$ outside imply $\beta<\gamma$ and $\delta>\epsilon$. 

If $BC$ and $DE$ are disjoint, then either both $B,C$ are outside $\triangle ABC$, or both $B,C$ are inside $\triangle ABC$. The case both $B,C$ are outside $\triangle ABC$ is the second of Figure \ref{geom2}. Since $\triangle ABC$ is an isosceles triangle, we get $\beta-\beta'=\gamma-\gamma'$. Therefore $\beta>\gamma$ if and only if $\beta'>\gamma'$. The case both $B,C$ are inside $\triangle ABC$ is the third of Figure \ref{geom2}. Since $\triangle ABC$ is an isosceles triangle, we get $\beta+\beta'=\gamma+\gamma'$. Therefore $\beta>\gamma$ if and only if $\beta'<\gamma'$. Since $\delta+\delta'=\gamma+\gamma'$, we also have $\delta<\epsilon$ if and only if $\delta'<\epsilon'$. In either case, the proof of Lemma \ref{geometry1} is reduced to the proof of the following similar result for quadrilaterals.

\begin{figure}[htp]
\centering
\begin{tikzpicture}


\draw
	(-1.2,0.6) 
	-- node[fill=white,inner sep=1] {\small $a$}
	(-0.8,-0.6)
	-- (0.8,-0.6)
	-- node[fill=white,inner sep=1] {\small $a$}
	(1.2,0.6)
	-- cycle;

\node at (-1.4,0.6) {\small $B$};
\node at (1.4,0.6) {\small $C$};
\node at (-1,-0.6) {\small $D$};
\node at (1,-0.6) {\small $E$};

\node at (-0.9,0.4) {\small $\beta$};
\node at (0.95,0.4) {\small $\gamma$};
\node at (-0.7,-0.4) {\small $\delta$};
\node at (0.7,-0.4) {\small $\epsilon$};

\end{tikzpicture}
\caption{Geometrical constraint for quadrilateral.}
\label{geom3}
\end{figure}
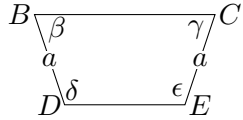

\begin{lemma}\label{geometry2}
If the spherical quadrilateral in Figure \ref{geom3} is simple and has a pair of equal edges $a$, then $\beta>\gamma$ is equivalent to $\delta<\epsilon$.
\end{lemma}

Similar to pentagon, we note that the reformulation of Lemma \ref{geometry2} in terms of the outside is equivalent to the original formulation in terms of the inside.

To prove Lemma \ref{geometry2}, we use the conformally accurate way of drawing great circles on the sphere. Let the circle $\Gamma$ be the stereographic projection (from the north pole to the tangent space of the south pole) of the equator. The antipodal points on the equator are then projected to the antipodal points on $\Gamma$. We denote the antipodal point of $P$ by $P^*$. Since the intersection of any great arc with the equator is a pair of antipodal points on the equator, the great circles of the sphere are in one-to-one correspondence with the circles (and straight lines) on the plane that intersect $\Gamma$ at a pair of antipodal points.

\begin{proof}
Suppose $a>\pi$. In Figure \ref{geom4}, we draw great circles ($\bigcirc BPP^*$ and $\bigcirc CPP^*$) containing the two $a$-edges. They intersect at a pair of antipodal points $P,P^*$ and divide the sphere into four $2$-gons. Since $a>\pi$ and the boundary of the quadrilateral is simple, $P,P^*$ lie in different $a$-edges. Up to symmetry, therefore, there are two ways the four vertices $B,C,D,E$ of the two $a$-edges can be located, described by the two pictures in Figure \ref{geom4}. Moreover, since $a>\pi$, the antipodal point $B^*$ of $B$ lies in the $a$-edge $BP^*D$.

In the first of Figure \ref{geom4}, we consider two great arcs connecting $B$ and $C$. One great arc (the solid one) is completely contained in the indicated $2$-gon. The other great arc (the dashed one) intersects the $a$-edge $BP^*D$ at the antipodal point $B^*$ and therefore cannot be an edge of the quadrilateral. We conclude that the $BC$ edge is the solid one. By the same reason, the edge $DE$ is also the solid one, that is completely contained in the indicated $2$-gon. Then the picture implies
\[
\beta<\pi<\gamma,\quad
\delta>\pi>\epsilon.
\]
Similar argument gives the $BC$ edge and $DE$ edges of the quadrilateral in the second of Figure \ref{geom4}, and we get the same inequalities above.

In all the subsequent argument, we may assume $a<\pi$.

\begin{figure}[htp]
\centering
\begin{tikzpicture}[scale=0.9]


\draw[dotted]
	(0,0) circle (1);
	
\coordinate (P) at (0,1);
\coordinate (Q) at (0,-1);
\coordinate (P1) at (160:1);
\coordinate (Q1) at (-20:1);
\coordinate (P2) at (180:1);
\coordinate (Q2) at (0:1);

\node at (0.05,1.25) {\small $P$};	
\node at (0,-1.3) {\small $P^*$};
\node at (45:1.2) {\small $\Gamma$};	


\begin{scope}[xshift=-1cm]

\coordinate (B) at (100:1.414);
\coordinate (D) at (10:1.414);
\coordinate (BB) at (-30:1.414);

\node at (105:1.6) {\small $B$};
\node at (110:1.2) {\small $\beta$};
\node at (0:1.6) {\small $D$};
\node at (20:1.3) {\small $\delta$};
\node at (-20:1.65) {\small $B^*$};

\end{scope}


\begin{scope}[xshift=1cm]

\coordinate (C) at (186.3:1.414);
\coordinate (E) at (-54.1:1.414);

\node at (-50:1.6) {\small $E$};
\node at (-45:1.25) {\small $\epsilon$};

\end{scope}


\begin{scope}
    \clip (-2.5,-1.8) rectangle (2.5,1.8);

\arcThroughThreePoints{B}{BB}{D};
\arcThroughThreePoints[dashed]{D}{P}{B};

\arcThroughThreePoints{E}{P}{C};
\arcThroughThreePoints[dashed]{C}{Q}{E};

\arcThroughThreePoints{B}{P1}{C};
\arcThroughThreePoints[dashed]{C}{Q1}{B};

\arcThroughThreePoints{E}{Q2}{D};
\arcThroughThreePoints[dashed]{D}{Q2}{E};

\end{scope}

\node[fill=white,inner sep=2] at (-2.4,0) {\small $a$};
\node[fill=white,inner sep=2] at (2.4,0) {\small $a$};

\node[fill=white,inner sep=0] at (-0.57,0.25) {\small $C$};
\node[fill=white,inner sep=1] at (-0.35,-0.35) {\small $\gamma$};


\begin{scope}[xshift=5.5cm]

\draw[dotted]
	(0,0) circle (1);
	
\coordinate (P) at (0,1);
\coordinate (Q) at (0,-1);
\coordinate (P1) at (135:1);
\coordinate (Q1) at (-45:1);
\coordinate (P2) at (180:1);
\coordinate (Q2) at (0:1);
	

\begin{scope}[xshift=-1cm]

\coordinate (B) at (100:1.414);
\coordinate (D) at (10:1.414);
\coordinate (BB) at (-30:1.414);

\node at (95:1.6) {\small $B$};
\node at (110:1.55) {\small $\beta$};
\node at (5:1.2) {\small $D$};
\node at (15:1.6) {\small $\delta$};

\end{scope}


\begin{scope}[xshift=1cm]

\coordinate (C) at (170:1.414);
\coordinate (EE) at (-54.1:1.414);
\coordinate (E) at (-69.4:1.414);

\node at (160:1.3) {\small $\epsilon$};
\node at (180:1.57) {\small $E$};

\end{scope}


\begin{scope}
    \clip (-2.5,-1.8) rectangle (2.5,1.8);

\arcThroughThreePoints{B}{BB}{D};
\arcThroughThreePoints[dashed]{D}{P}{B};

\arcThroughThreePoints{E}{P}{C};
\arcThroughThreePoints[dashed]{C}{Q}{E};

\arcThroughThreePoints{E}{P1}{B};
\arcThroughThreePoints[dashed]{B}{P1}{E};

\arcThroughThreePoints{D}{EE}{C};
\arcThroughThreePoints[dashed]{C}{EE}{D};

\end{scope}

\node[fill=white,inner sep=2] at (-2.4,0) {\small $a$};
\node[fill=white,inner sep=2] at (2.4,0) {\small $a$};

\node[fill=white,inner sep=0] at (1.9,-1.33) {\small $C$};
\node[fill=white,inner sep=1] at (1.3,-1.3) {\small $\gamma$};

\end{scope}


\begin{scope}[xshift=10cm]

\draw
	(0,0) circle (1.4);

\coordinate (D) at (0,1.4);
\coordinate (DD) at (0,-1.4);
\coordinate (E) at (-60:1.4);
\coordinate (EE) at (120:1.4);

\coordinate (B) at (50:0.8);

\coordinate (C) at (200:0.52);

\coordinate (P) at (220:1.4);
\coordinate (Q) at (40:1.4);


\begin{scope}
    \clip (-1.5,-1.8) rectangle (1.5,1.8);
    
\arcThroughThreePoints{B}{DD}{D};
\arcThroughThreePoints[dashed]{D}{DD}{B};

\arcThroughThreePoints{B}{P}{C};
\arcThroughThreePoints[dashed]{C}{P}{B};

\arcThroughThreePoints{C}{EE}{E};
\arcThroughThreePoints[dashed]{E}{EE}{C};

\end{scope}

\node at (0.1,1.6) {\small $D$};
\node at (0.1,-1.6) {\small $D^*$};
\node at (93:1.2) {\small $\delta$};

\node[fill=white,inner sep=1] at (40:0.95) {\small $B$};
\node at (65:0.7) {\small $\beta$};

\node at (215:0.3) {\small $C$};
\node[fill=white,inner sep=1] at (195:0.68) {\small $\gamma$};

\node at (-55:1.2) {\small $E$};
\node at (-75:1.25) {\small $\epsilon$};

\end{scope}
           
\end{tikzpicture}
\caption{Lemma \ref{geometry2}: $a>\pi$, and $DE<\pi$ in case $\delta,\gamma<\pi$.}
\label{geom4}
\end{figure}
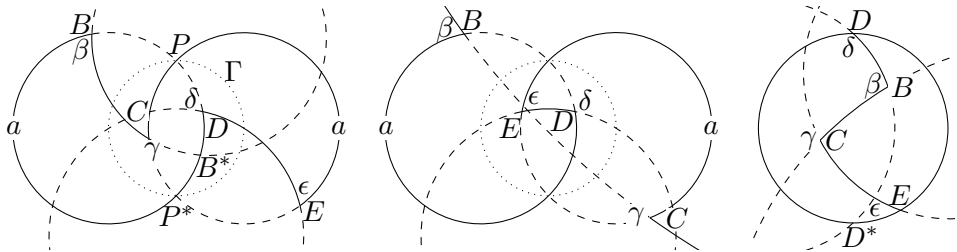

Next we argue that, if $\delta,\epsilon<\pi$, then we may further assume $DE<\pi$. Suppose $DE>\pi$. We draw the great circle $\bigcirc DE$ containing $DE$ in the third of Figure \ref{geom4}. Since $\delta,\epsilon<\pi$, both $DB$ and $EC$ lie in the same hemisphere $H$ bounded by the circle. This implies that both $B$ and $C$ lie in $H$. Therefore among two great arcs connecting $B$ and $C$, one intersects $\bigcirc DE$ at two antipodal points. By $DE>\pi$, one such point lies on $DE$, and therefore this arc is not the $BC$ edge. This proves that the $BC$ edge also lies in the hemisphere $H$. Since all three edges connecting $D,B,C,E$ lie in $H$, we can have the complement of the quadrilateral in $H$. The complement is still a quadrilateral, which differs from the original one by replacing $DE$ with the other great arc (of length $2\pi-DE<\pi$) connecting $D$ and $E$, and by replacing $\beta,\gamma,\delta,\epsilon$ with $2\pi-\beta,2\pi-\gamma,\pi-\delta,\pi-\epsilon$. Moreover, we still have $\pi-\delta,\pi-\epsilon<\pi$. It is easy to see that the lemma for the new quadrilateral is equivalent to the lemma for the original quadrilateral.

Now we are ready to prove the lemma. We divide into two cases. The first is at least three angles $<\pi$. By considering the complement quadrilateral, this also covers the case of at least three angles $>\pi$. The second is the case of two angles $>\pi$ and two angles $<\pi$. 

Suppose at least three angles $<\pi$. Up to symmetry, we may assume $\beta,\delta,\epsilon<\pi$. By the earlier argument, we may further assume that $a<\pi$ and $DE<\pi$. In Figure \ref{geom5}, we draw the great circles $\bigcirc DB$ and $\bigcirc DE$ containing $DB$ and $DE$. The two great circles divide the sphere into four $2$-gons. By $\delta<\pi$, we may assume that $\delta$ is an angle of the middle $2$-gon. By $DB=a<\pi$ and $DE<\pi$, $B$ and $E$ lie on the two edges of the middle $2$-gon. By $\beta,\epsilon<\pi$, we find that $BC$ and $EC$ are inside the middle $2$-gon. We also prolong the $a$-edge $EC$ to intersect the boundary of middle $2$-gon at $T$. The two pictures in Figure \ref{geom5} refer to $\gamma>\pi$ and $\gamma<\pi$. In the first picture, we have $DT<DB=a=EC<ET$. Since all angles of $\triangle DET$ are $<\pi$, this implies that $\epsilon=\angle DET<\angle EDT=\delta$. We also have $\beta<\pi<\gamma$. In the second picture, the great circles $\bigcirc DB$ and $\bigcirc EC$ containing the two $a$-edges intersect at antipodal points $T$ and $T^*$, and the two antipodal points do not lie on the two $a$-edges. This implies $BT+a+DT^*=\pi=CT+a+ET^*$. Since all angles in $\triangle BCT$ and $\triangle DET^*$ are $<\pi$, we then have 
\[
\beta>\gamma
\iff a+ET^*>a+DT^*
\iff a+CT<a+BT
\iff \delta<\epsilon.
\]

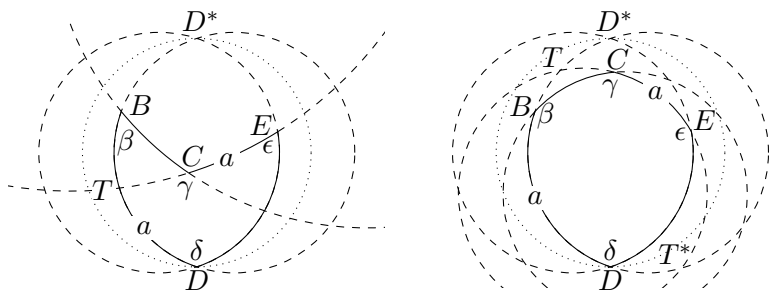
\begin{figure}[htp]
\centering
\begin{tikzpicture}


\draw[dashed]
	(-0.5,0) circle (1.6)
	(0.5,0) circle (1.6);
	
\draw[dotted]
	(0,0) circle (1.52);
			
\coordinate (DD) at (0,1.52);
\coordinate (D) at (0,-1.52);
\coordinate (C) at (-0.108,-0.27);
\coordinate (P1) at (145:1.52);
\coordinate (Q1) at (-35:1.52);
\coordinate (P2) at (200:1.52);
\coordinate (Q2) at (20:1.52);

\node at (0.05,1.75) {\small $D^*$};	
\node at (0,-1.7) {\small $D$};	
\node at (0,-1.3) {\small $\delta$};


\begin{scope}[xshift=-0.5cm]

\coordinate (E) at (10:1.6);

\node at (15:1.4) {\small $E$};
\node at (3:1.45) {\small $\epsilon$};

\end{scope}


\begin{scope}[xshift=0.5cm]

\coordinate (B) at (160:1.6);

\node at (153:1.4) {\small $B$};
\node at (175:1.45) {\small $\beta$};

\end{scope}


\begin{scope}
    \clip (-2.5,-1.8) rectangle (2.5,1.8);

\arcThroughThreePoints{D}{DD}{E};
\arcThroughThreePoints{B}{DD}{D};

\arcThroughThreePoints{B}{P1}{C};
\arcThroughThreePoints[dashed]{C}{P1}{B};

\arcThroughThreePoints{C}{Q2}{E};
\arcThroughThreePoints[dashed]{E}{Q2}{C};

\end{scope}

\node at (-0.05,-0.05) {\small $C$};
\node at (-0.15,-0.5) {\small $\gamma$};
\node[fill=white,inner sep=0] at (-1.25,-0.5) {\small $T$};

\node[fill=white,inner sep=2] at (0.4,-0.1) {\small $a$};
\node[fill=white,inner sep=2] at (-0.7,-1) {\small $a$};


\begin{scope}[xshift=5.5cm]

\draw[dashed]
	(-0.5,0) circle (1.6)
	(0.5,0) circle (1.6);
\draw[dotted]
	(0,0) circle (1.52);
			
\coordinate (DD) at (0,1.52);
\coordinate (D) at (0,-1.52);
\coordinate (C) at (0.05,1.07);
\coordinate (P1) at (145:1.52);
\coordinate (Q1) at (-35:1.52);
\coordinate (P2) at (200:1.52);
\coordinate (Q2) at (20:1.52);

\node at (0.05,1.75) {\small $D^*$};	
\node at (0,-1.7) {\small $D$};	
\node at (0,-1.3) {\small $\delta$};


\begin{scope}[xshift=-0.5cm]

\coordinate (E) at (10:1.6);

\node at (14:1.8) {\small $E$};
\node at (10:1.45) {\small $\epsilon$};

\end{scope}


\begin{scope}[xshift=0.5cm]

\coordinate (B) at (160:1.6);

\node at (160:1.8) {\small $B$};
\node at (160:1.45) {\small $\beta$};

\node at (135:1.8) {\small $T$};
\node at (-75:1.4) {\small $T^*$};

\end{scope}


\begin{scope}
    \clip (-2.5,-1.8) rectangle (2.5,1.8);

\arcThroughThreePoints{D}{DD}{E};
\arcThroughThreePoints{B}{DD}{D};

\arcThroughThreePoints{C}{P2}{B};
\arcThroughThreePoints[dashed]{B}{P2}{C};

\arcThroughThreePoints{E}{P1}{C};
\arcThroughThreePoints[dashed]{C}{P1}{E};

\end{scope}

\node at (0.1,1.25) {\small $C$};
\node at (0,0.85) {\small $\gamma$};

\node[fill=white,inner sep=2] at (0.6,0.8) {\small $a$};
\node[fill=white,inner sep=2] at (-1,-0.6) {\small $a$};

\end{scope}
           
\end{tikzpicture}
\caption{Lemma \ref{geometry2}: At least three angles $>\pi$.}
\label{geom5}
\end{figure}

Suppose two angles $>\pi$ and two angles $<\pi$. Then up to symmetry, we need to consider the following three cases.
\begin{enumerate}
\item $\beta,\gamma>\pi$ and $\delta,\epsilon<\pi$.
\item $\beta,\delta>\pi$ and $\gamma,\epsilon<\pi$.
\item $\beta,\epsilon>\pi$ and $\gamma,\delta<\pi$.
\end{enumerate}

In the first case, by the earlier argument, we may additionally assume $a<\pi$ and $DE<\pi$. We draw great circles $\bigcirc BD$ and $\bigcirc CE$ containing the two $a$-edges. Since $a<\pi$ and the two $a$-edges do not intersect, we may assume that they are inside the two edges of the middle $2$-gon bounded by $\bigcirc BD$ and $\bigcirc CE$, as in the first of Figure \ref{geom6}. Since $DE<\pi$, the edge $DE$ lies in the middle $2$-gon. Since $\delta,\epsilon<\pi$, the two angles also lies in the middle $2$-gon. Then $\beta,\gamma>\pi$ implies that the $BC$ edge lies outside the middle $2$-gon. If we replace the $BC$ edge by the other great arc $c$ connecting the two points, then we get a new quadrilateral with angles $\beta-\pi,\gamma-\pi,\delta,\epsilon$. All four angles of the new quadrilateral are $<\pi$, and the earlier argument shows that the lemma holds for the new quadrilateral. Then it is easy to see that the lemma for the new quadrilateral is equivalent to the lemma for the original quadrilateral. This completes the proof of the case.

\begin{figure}[htp]
\centering
\begin{tikzpicture}


\draw[dashed]
	(-0.5,0) circle (1.6)
	(0.5,0) circle (1.6);

\fill 
	(0.96,-0.65) circle (0.03);
	
\coordinate (AA) at (0,1.52);
\coordinate (A) at (0,-1.52);
\coordinate (E) at (0.96,-0.65);
\coordinate (B) at (-0.96,0.65);
\coordinate (P1) at (145:1.52);
\coordinate (Q1) at (-35:1.52);
\coordinate (P2) at (200:1.52);
\coordinate (Q2) at (20:1.52);



\begin{scope}[xshift=-0.5cm]

\coordinate (C) at (40:1.6);

\node at (36:1.4) {\small $C$};
\node at (-25:1.8) {\small $E$};
\node at (-22:1.45) {\small $\epsilon$};

\end{scope}


\begin{scope}[xshift=0.5cm]

\coordinate (D) at (220:1.6);

\node at (160:1.45) {\small $B$};
\node at (220:1.8) {\small $D$};
\node at (217:1.4) {\small $\delta$};

\end{scope}


\begin{scope}
    \clip (-2.5,-1.8) rectangle (2.5,1.8);

\arcThroughThreePoints{E}{AA}{C};
\arcThroughThreePoints{B}{AA}{D};

\arcThroughThreePoints{D}{P2}{E};
\arcThroughThreePoints[dashed]{E}{P2}{D};

\arcThroughThreePoints{B}{Q2}{C};
\arcThroughThreePoints[dashed]{C}{Q2}{B};

\end{scope}

\node[fill=white,inner sep=0] at (-0.85,0.85) {\small $\beta$};
\node[fill=white,inner sep=0] at (0.6,1.1) {\small $\gamma$};

\node[fill=white,inner sep=2] at (1.05,0.3) {\small $a$};
\node[fill=white,inner sep=2] at (-1.05,-0.3) {\small $a$};

\node[fill=white,inner sep=2] at (-0.1,1.1) {\small $c$};


\begin{scope}[xshift=5cm]

\draw
	(0,0) circle (1.6);
	
\coordinate (EE) at (0,1.6);
\coordinate (E) at (0,-1.6);
\coordinate (CC) at (200:1.6);
\coordinate (C) at (20:1.6);
\coordinate (B) at (0.7,0);
\coordinate (D) at (-0.31,0.98);
\coordinate (P) at (130:1.6);
\coordinate (Q) at (-50:1.6);		


\begin{scope}
    \clip (-1.8,-1.8) rectangle (1.8,1.8);

\arcThroughThreePoints{B}{CC}{C};
\arcThroughThreePoints[dashed]{C}{CC}{B};

\arcThroughThreePoints{D}{EE}{E};
\arcThroughThreePoints[dashed]{E}{EE}{D};

\arcThroughThreePoints{B}{P}{D};
\arcThroughThreePoints[dashed]{D}{P}{B};

\end{scope}

\node at (92:1.8) {\small $E^*$};	
\node at (-94:1.8) {\small $E$};	
\node at (16:1.8) {\small $C$};
\node at (205:1.85) {\small $C^*$};
\node at (0.8,0.3) {\small $B$};
\node at (-0.4,1.25) {\small $D$};

\node at (0.62,-0.25) {\small $\beta$};
\node at (10:1.45) {\small $\gamma$};
\node at (-0.25,0.7) {\small $\delta$};
\node at (-87:1.45) {\small $\epsilon$};

\node[fill=white,inner sep=2] at (-30:1.6) {\small $a$};
\node[fill=white,inner sep=1] at (0.3,0.5) {\small $a$};

\end{scope}
         
\end{tikzpicture}
\caption{Lemma \ref{geometry2}: Two angles $<\pi$ and two angles $>\pi$.}
\label{geom6}
\end{figure}
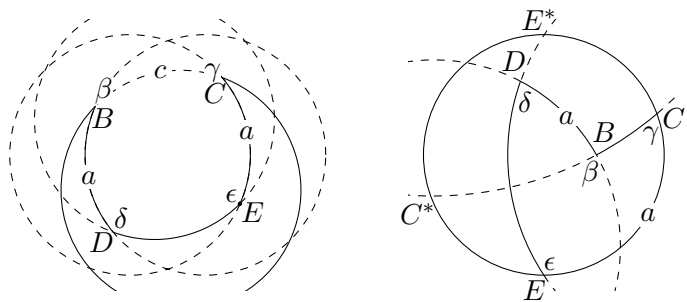

In the second case, we may again assume $a<\pi$. We draw the circle $\bigcirc CE$ containing the $a$-edge $CE$, as in the second of Figure \ref{geom6}. Since $\gamma,\epsilon<\pi$, the two angles lie in the same hemisphere $H_1$ bounded by $\bigcirc CE$. Then the prolongation of $CB$ and $ED$ intersect at a point inside the hemisphere $H_1$. Since $CB$ and $ED$ do not intersect, we must have either $B,C$ lie in the same hemisphere $H_2$ bounded by the circle $\bigcirc DEE^*$, or $D,E$ lie in the same hemisphere bounded by the circle $\bigcirc BCC^*$. Without loss of generality, we may assume the first scenario happens, as in the second of Figure \ref{geom6}. Since $a<\pi$, we also know that both $BC$ and $CE$ lie in $H_2$, as in the second of Figure \ref{geom6}. Now $BD=a<\pi$ implies that the $BD$ edge cannot intersect $\bigcirc CE$ at two points. Therefore the $BD$ edge also lies inside $H_1$. This implies that $\delta$ is the angle between an edge $BD$ inside $H_2$ and the boundary $\bigcirc DEE^*$ of $H_2$. Such angle is always $<\pi$. This contradicts the assumption that $\delta>\pi$.

Finally, the third case is consistent with the conclusion of lemma.
\end{proof}

\section{Spherical Trigonometry}
\label{constraint2}

\begin{figure}[htp]
\centering
\begin{tikzpicture}


\draw[dashed]
	(-1,0) -- node[fill=white,inner sep=1] {\small $x$}
	(1,0);

\node at (0,0.5) {\small $\alpha$};	
\node[fill=white,inner sep=1] at (-0.75,-0.05) {\small $\beta$};
\node[fill=white,inner sep=1] at (0.8,-0.1) {\small $\gamma$};	
\node at (-0.5,-0.8) {\small $\delta$};	
\node at (0.5,-0.8) {\small $\epsilon$};

\draw
	(0,0.7) -- node[fill=white,inner sep=1] {\small $b$}
	(-1,0) -- node[fill=white,inner sep=1] {\small $a$}
	(-0.6,-1) -- node[fill=white,inner sep=1] {\small $a$}
	(0.6,-1) -- node[fill=white,inner sep=1] {\small $a$}
	(1,0) -- node[fill=white,inner sep=1] {\small $b$}
	cycle;


\begin{scope}[xshift=3cm]

\draw
	(0,0.7) -- node[fill=white,inner sep=1] {\small $b$}
	(-1,-1) -- node[fill=white,inner sep=1] {\small $x$}
	(1,-1) -- node[fill=white,inner sep=1] {\small $b$}
	(0,0.7);

\node at (0,0.4) {\small $\alpha$};
\node at (-0.7,-0.8) {\small $\theta$};
\node at (0.7,-0.8) {\small $\theta$};

\end{scope}


\begin{scope}[xshift=6cm]

\draw
	(1.2,0.7) -- node[fill=white,inner sep=1] {\small $x$}
		(-1.2,0.7) -- node[fill=white,inner sep=1] {\small $a$}
		(-0.8,-1) -- node[fill=white,inner sep=1] {\small $a$}
		(0.8,-1) -- node[fill=white,inner sep=1] {\small $a$} 
		(1.2,0.7);

\node at (-0.65,0.5) {\small $\beta-\theta$};
\node at (0.7,0.5) {\small $\gamma-\theta$};
\node at (-0.7,-0.8) {\small $\delta$};
\node at (0.7,-0.8) {\small $\epsilon$};

\end{scope}


\begin{scope}[xshift=9cm]

\draw
	(1.2,0.7) -- node[fill=white,inner sep=1] {\small $x$}
	(-1.2,0.7) -- node[fill=white,inner sep=1] {\small $a$}
	(-0.8,-1) -- node[fill=white,inner sep=1] {\small $a$}
	(0.8,-1) -- node[fill=white,inner sep=1] {\small $a$} 
	(1.2,0.7);
		
\draw[dashed]
	(-0.8,-1) -- node[fill=white,inner sep=1] {\small $y$}
	(1.2,0.7);

\node at (-0.95,0.5) {\small $\beta$};
\node at (0.7,-0.8) {\small $\epsilon$};

\node at (0.9,0.2) {\small $\phi$};
\node at (-0.35,-0.8) {\small $\phi$};

\node[rotate=25] at (0.48,0.4) {\scriptsize $\gamma-\phi$};
\node[rotate=50] at (-0.5,-0.5) {\scriptsize $\delta-\phi$};

\node[draw,shape=circle, inner sep=1] at (0.4,-0.5) {\small 1};
\node[draw,shape=circle, inner sep=1] at (-0.4,0.2) {\small 2};

\end{scope}

\end{tikzpicture}
\caption{Spherical trigonometry of pentagon of $a^3b^2$ type.}
\label{geom7}
\end{figure}
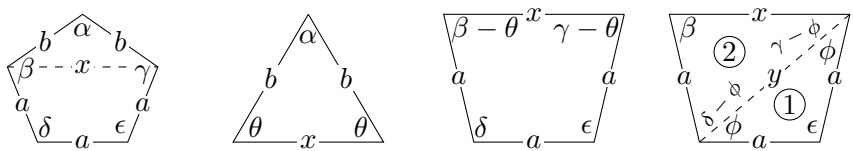

Lemma \ref{geometry1} is a very rough constraint on spherical pentagons because it is about inequalities. For example, for given values of $a,b$ and five angles satisfying the lemma, there is no guarantee that such pentagon exists. In this section, we give more precise equalities that must be satisfied by spherical pentagon in the first of Figure \ref{geom7} (i.e., the third type $a^3b^2$ in Figure \ref{edges1}). Note that the results of this section do not require $a,b$ to be distinct, and are therefore also applicable to the edge combination $a^5$.

A general spherical pentagon allows $7$ free parameters. The spherical pentagon in the first of Figure \ref{geom7} satisfies $3$ independent edge length equalities. Therefore the pentagon allows $7-3=4$ free parameters. This means that four angles ($\beta,\gamma,\delta,\epsilon$, for example) completely determine the pentagon.

We divide the pentagon into an isosceles triangle (second of Figure \ref{geom7}) and a quadrilateral (third of Figure \ref{geom7}). If $\alpha<\pi$, then Figure \ref{geom2} gives three possible scenarios the pentagon is divided. If $\alpha>\pi$, then Figure \ref{geom2} gives three possible divisions of the outside pentagon. The subsequent discussion is guided by the first of Figure \ref{geom7} (the second of Figure \ref{geom2}), and further verification shows that the conclusions are still valid in all cases.

A general triangle allows $3$ free parameters. The isosceles triangle introduces $1$ equality and allows $3-1=2$ free parameters. Therefore $\alpha$ and $\theta$ completely determine the isosceles triangle. The edge length $x$ is given by
\[
\cos x=\frac{\cos\alpha+\cos^2\theta}{\sin^2\theta}
=\cos\alpha+(1+\cos\alpha)\cot^2\theta.
\]
The edge length $b$ is further given by
\[
\sin x\cos\theta=\sin b\cos b(1-\cos\alpha),\quad
\sin x\sin\theta=\sin b\sin\alpha.
\]

For the quadrilateral in the third of Figure \ref{geom7}, we relabel $\beta-\theta,\gamma-\theta$ by $\beta,\gamma$ and get the fourth of Figure \ref{geom7}. A general quadrilateral allows $5$ free parameters. Three equal edges introduces $2$ equalities and leaves $5-2=3$ free parameters. Therefore its four angles satisfy one equality, and the edge lengths $a$ and $x$ can be expressed in terms of the four angles.

We add an arc $y$ to the quadrilateral and get two triangles. By the sine law for the triangle labeled 1 and by applying the known formulae to the isosceles triangle labeled 2, we have
\begin{align*}
\sin a\sin\beta
&=\sin y\sin(\gamma-\phi) \\
&=\sin y\sin\gamma\cos\phi-\sin y\cos\gamma\sin\phi \\
&=\sin a\cos a\sin\gamma(1-\cos\epsilon)-\sin a\cos\gamma\sin\epsilon.
\end{align*}
Canceling $\sin a$, we get
\begin{equation}\label{eq3}
\sin\beta+\cos\gamma\sin\epsilon=\cos a\sin\gamma(1-\cos\epsilon).
\end{equation}
By the symmetry of exchanging $\beta$ and $\gamma$, and exchanging $\delta$ and $\epsilon$, we get
\begin{equation}\label{eq4}
\sin\gamma+\cos\beta\sin\delta=\cos a\sin\beta(1-\cos\delta).
\end{equation}
Canceling $\cos a$ from \eqref{eq3} and \eqref{eq4}, we get an equality for the four angles of the quadrilateral
\[
(\sin\beta+\cos\gamma\sin\epsilon)\sin\beta(1-\cos\delta)
=(\sin\gamma+\cos\beta\sin\delta)\sin\gamma(1-\cos\epsilon).
\]
Then $a$ is determined by the four angles through \eqref{eq3} or \eqref{eq4}. Moreover, $x$ is determined below
\begin{align}
\cos x
&=\cos a\cos y+\sin a\sin y\cos(\delta-\phi) \nonumber \\
&=\cos a\cos y+\sin a\cos\delta\sin y\cos \phi+\sin a\sin\delta\sin y\sin \phi \nonumber \\
&=\cos a(\cos^2a+\sin^2a\cos\epsilon) \nonumber \\
&\quad +\sin a\cos\delta\sin a\cos a(1-\cos\epsilon)+\sin a\sin\delta\sin a\sin\epsilon \nonumber \\
&=\cos a-\sin^2a\cos a(1-\cos\delta)(1-\cos\epsilon)+\sin^2a\sin\delta\sin\epsilon \nonumber \\
&=\cos^3a(1-\cos\delta)(1-\cos\epsilon)
-\cos^2a\sin\delta\sin\epsilon \nonumber \\
&\quad +\cos a(\cos\delta+\cos\epsilon-\cos\delta\cos\epsilon)
+\sin\delta\sin\epsilon. \label{eq6}
\end{align}

Back to the pentagon in the first of Figure \ref{geom7}, we substitute $\beta,\gamma,\delta,\epsilon$ in \eqref{eq3} and \eqref{eq4} by $\beta-\theta,\gamma-\theta,\delta,\epsilon$. After dividing by $\sin\theta$. We get
\begin{align*}
&(\sin\beta+\cos\gamma\sin\epsilon)\cot\theta-\cos\beta+\sin\gamma\sin\epsilon \nonumber \\
&=\cos a(\sin\gamma\cot\theta-\cos\gamma)(1-\cos\epsilon),  \\
&(\sin\gamma+\cos\beta\sin\delta)\cot\theta-\cos\gamma+\sin\beta\sin\delta \nonumber \\
&=\cos a(\sin\beta\cot\theta-\cos\beta)(1-\cos\delta). 
\end{align*}
This is a system of two equations for $\cos a$ and $\cot\theta$. The solutions are quadratic equation for $\cos a$ and $\cot\theta$
\begin{align}
L\cos^2a+M\cos a+N &=0, \label{eq10} \\
P\cot^2\theta+Q\cot\theta+R &=0. \label{eq11}
\end{align}
The coefficients are functions of $\beta,\gamma,\delta,\epsilon$
\begin{align*}
L
&=\sin(\beta-\gamma)(1-\cos\delta)(1-\cos\epsilon), \\
M
&=\cos(\beta-\gamma)(\sin\epsilon-\sin\delta+\sin(\delta-\epsilon)),  \\
N
&=\sin\delta-\sin\epsilon-\sin(\beta-\gamma)(1-\sin\delta\sin\epsilon), \\
P
&=\sin\beta(\sin\beta+\cos\gamma\sin\epsilon)(1-\cos\delta) \\
&\quad -\sin\gamma(\sin\gamma+\cos\beta\sin\delta)(1-\cos\epsilon), \\
Q
&=-(\sin 2\beta+\cos(\beta+\gamma)\sin\epsilon)(1-\cos\delta) \\
&\quad  +(\sin 2\gamma+\cos(\beta+\gamma)\sin\delta)(1-\cos\epsilon),  \\
R
&=\cos\beta(\cos\beta-\sin\gamma\sin\epsilon)(1-\cos\delta) \\
&\quad -\cos\gamma(\cos\gamma-\sin\beta\sin\delta)(1-\cos\epsilon).
\end{align*}
Combined with \eqref{eq6} that determines $x$, we may further determine $\alpha$ and $b$ of the isosceles triangle.

\section{Neighborhood Tiling}
\label{nd}

Our strategy for constructing tilings is to first construct tilings of the neighborhood of the special tile in Proposition \ref{base_tile}. In this paper, we concentrate on the neighborhood of a $3^5$-tile. The neighborhood tiling implies some useful information about edges and angles. The information is already enough for dismissing the edge combinations $a^2b^2c$ and $a^3bc$. 

\begin{theorem}\label{thm1}
If an edge-to-edge spherical tiling by congruent pentagons has a $3^5$-tile and the pentagon has edge length combination $a^2b^2c$, with $a,b,c$ distinct, then the number of tiles is $f=12$.
\end{theorem} 

\begin{figure}[htp]
\centering
\begin{tikzpicture}[>=latex]


\foreach \a in {0,1,2}
{
\begin{scope}[xshift=3.5*\a cm]

\foreach \x in {1,...,5}
{
\coordinate (A\x X\a) at (-54+72*\x:0.7);
\coordinate (B\x X\a) at (-54+72*\x:1.3);
\coordinate (C\x X\a) at (-18+72*\x:1.7);
\coordinate (P\x X\a) at (-18+72*\x:1.7);
}

\node[draw,shape=circle, inner sep=1] at (0,0) {\small 1};
\foreach \x in {2,...,6}
\node[draw,shape=circle, inner sep=1] at (-90+72*\x:1) {\small \x};

\end{scope}
}

\draw
	(A2X0) -- (A3X0) -- (A4X0) 
	(A2X0) -- (B2X0) -- (C1X0)
	(A5X0) -- (B5X0) -- (C5X0)
	(B5X0) -- (C4X0)
	(A3X0) -- (B3X0)
	(A2X1) -- (A3X1) -- (A4X1) 
	(B1X1) -- (C1X1) -- (B2X1)
	(A3X1) -- (B3X1)
	(A3X2) -- (A4X2) -- (A5X2) -- (A1X2) 
	(A2X2) -- (B2X2)
	(B1X2) -- (C1X2) -- (B2X2) -- (C2X2) -- (B3X2)
	(A4X2) -- (B4X2) -- (C3X2)
	(A5X2) -- (B5X2) -- (C5X2)
	;

\draw[line width=1.5]
	(A2X0) -- (A1X0) -- (A5X0) 
	(A1X0) -- (B1X0)
	(A4X0) -- (B4X0)
	(C3X0) -- (B4X0) -- (C4X0)
	(B2X1) -- (C2X1) -- (B3X1) 
	(A2X1) -- (A1X1) -- (A5X1)
	(A1X1) -- (B1X1)
	(A3X2) -- (A2X2) 
	(A1X2) -- (B1X2)
	(B3X2) -- (C3X2)
	;

\draw[dashed]
	(A4X0) -- (A5X0) 
	(C1X0) -- (B1X0) -- (C5X0)
	(B3X0) -- (C3X0)
	(A2X1) -- (B2X1) 
	(A4X1) -- (A5X1)
	(A1X2) -- (A2X2)  
	(A3X2) -- (B3X2)
	(C5X2) -- (B1X2)
	; 
	
\draw[dotted]
	(B2X0) -- (C2X0) -- (B3X0) 
	(B3X1) -- (C3X1) -- (B4X1) -- (C4X1) -- (B5X1) -- (C5X1) -- (B1X1) 
	(A4X1) -- (B4X1)
	(A5X1) -- (B5X1)
	(B4X2) -- (C4X2) -- (B5X2) 
	;

\begin{scope}[xshift=3.5cm]

\node at (90:0.45) {\small $\alpha$}; 
\node at (162:0.45) {\small $\beta$};
\node at (18:0.45) {\small $\gamma$};
\node at (234:0.45) {\small $\delta$};
\node at (-54:0.45) {\small $\epsilon$};
\node at (78:0.8) {\small $\epsilon$}; 
\node at (32:0.8) {\small $\gamma$};
\node at (28:1.15) {\small $\alpha$};
\node at (82:1.15) {\small $\delta$};
\node at (54:1.45) {\small $\beta$};
\node at (102:0.8) {\small $\delta$}; 
\node at (148:0.8) {\small $\beta$};
\node at (152:1.15) {\small $\alpha$};
\node at (97:1.15) {\small $\epsilon$};
\node at (126:1.45) {\small $\gamma$};
\node at (176:0.8) {\small $\beta$}; 
\node at (4:0.8) {\small $\gamma$}; 

\end{scope}


\begin{scope}[xshift=7cm]

\node at (90:0.45) {\small $\alpha$}; 
\node at (162:0.45) {\small $\beta$};
\node at (18:0.45) {\small $\gamma$};
\node at (234:0.45) {\small $\delta$};
\node at (-54:0.45) {\small $\epsilon$};
\node at (78:0.8) {\small $\gamma$}; 
\node at (30:0.8) {\small $\alpha$};
\node at (28:1.15) {\small $\beta$};
\node at (82:1.15) {\small $\epsilon$};
\node at (54:1.45) {\small $\delta$};
\node at (102:0.8) {\small $\beta$}; 
\node at (148:0.8) {\small $\alpha$};
\node at (152:1.15) {\small $\gamma$};
\node at (97:1.15) {\small $\delta$};
\node at (126:1.5) {\small $\epsilon$};
\node at (174:0.8) {\small $\gamma$}; 
\node at (220:0.75) {\small $\epsilon$};
\node at (224:1.15) {\small $\delta$};
\node at (196:1.45) {\small $\beta$};
\node at (170:1.15) {\small $\alpha$};
\node at (2:0.75) {\small $\beta$}; 
\node at (-40:0.75) {\small $\delta$};
\node at (-44:1.15) {\small $\epsilon$};
\node at (-18:1.45) {\small $\gamma$};
\node at (10:1.15) {\small $\alpha$};

\end{scope}

\end{tikzpicture}
\caption{Neighborhood tilings for $a^2b^2c$ and $a^3bc$.}
\label{2a2bc_3abc}
\end{figure}
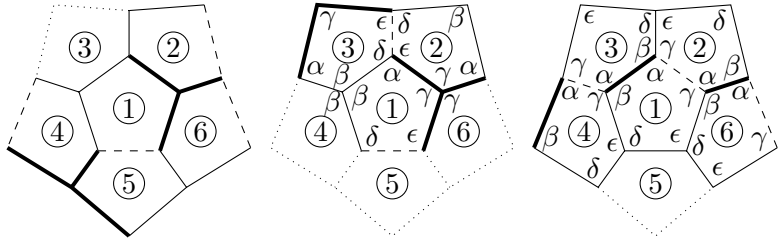

\begin{proof}
By Propositions \ref{edge_combo} and \ref{edge_combo2} (also see \cite[Proposition 7]{gsy}), the edge lengths of the pentagon are arranged as the second of Figure \ref{edges1}. A $3^5$-tile has another $5$ tiles around it, and its neighborhood is combinatorially given by Figure \ref{2a2bc_3abc}. We denote the tile by $P_1$ and its five neighboring tiles by $P_2,\dots,P_6$. 

The pentagon with edge combination $a^2b^2c$ is the last in Figure \ref{edges1}. Up to the combinatorial symmetry of the neighborhood, we may assume that the edges of $P_1$ are given as indicated, and the edge shared by $P_2,P_3$ is $a$ or $c$. The first of Figure \ref{2a2bc_3abc} is the case that the edge shared by $P_2,P_3$ is $a$. By the edge length consideration, we may successively determine all edges of $P_2,P_6,P_5,P_4$. Then we find three $a$-edges in $P_3$, a contradiction. The second of Figure \ref{2a2bc_3abc} is the case that the edge shared by $P_2,P_3$ is $c$. We denote the five angles of the pentagon by $\alpha$ ($ab$-angle), $\beta$ ($a^2$-angle), $\gamma$ ($b^2$-angle), $\delta$ ($ac$-angle), $\epsilon$ ($bc$-angle). By edge length consideration, we may determine (all edges and angles of) $P_2,P_3$. Then we may further determine the $a^2$-angle $\beta$ of $P_4$ and $b^2$-angle $\gamma$ of $P_6$. The angle sums at the three vertices imply
\[
\alpha+\delta+\epsilon=3\beta=3\gamma=2\pi.
\]
This further implies $\alpha+\beta+\gamma+\delta+\epsilon=\frac{10}{3}\pi$. By the angle sum equation \eqref{asump} for pentagon, we conclude $f=12$. 
\end{proof}

\begin{theorem}\label{thm2}
If an edge-to-edge spherical tiling by congruent pentagons has a $3^5$-tile and the pentagon has edge length combination $a^3bc$, with $a,b,c$ distinct, then the number of tiles is $f=12$.
\end{theorem}

\begin{proof}
By Proposition \ref{edge_combo} (also see \cite[Proposition 7]{gsy}), the edge lengths of the pentagon are arranged as the third of Figure \ref{edges1}. Up to symmetry, we may assume that the edges of $P_1$ is given by the third of Figure \ref{2a2bc_3abc}, and with all angles as indicated. Since each tile has only one $b$-edge and one $c$-edge, the edge shared by $P_2,P_3$ must be $a$. This determines (all edges and angles of) $P_2,P_3$. Then we may further determine $P_4,P_6$. 

Given the three known $a$-edges of $P_5$, there are two ways of arranging (edges and angles of) $P_5$. Either way gives two vertices $\delta^2\epsilon$ and $\delta\epsilon^2$ shared by $P_1,P_4,P_5$ and by $P_1,P_5,P_6$. The angle sums at the two vertices and at the vertex $\alpha\beta\gamma$ imply $\delta=\epsilon=\frac{2}{3}\pi$ and $\alpha+\beta+\gamma=2\pi$. Then we get $\alpha+\beta+\gamma+\delta+\epsilon=\frac{10}{3}\pi$. By the angle sum equation \eqref{asump} for pentagon, we conclude $f=12$.  
\end{proof}

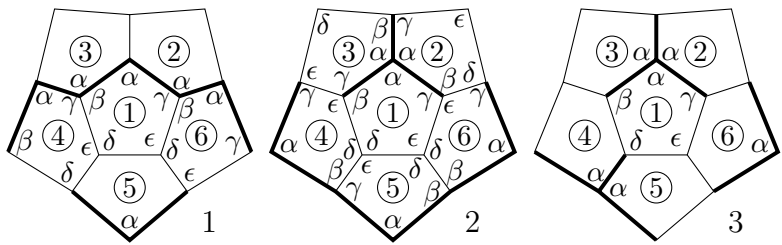
\begin{figure}[htp]
\centering
\begin{tikzpicture}[>=latex,scale=1]

\foreach \a in {0,1,2}
{
\begin{scope}[xshift=3.5*\a cm]

\foreach \x in {1,...,5}
\draw[rotate=72*\x]
	(90:0.7) -- (18:0.7) -- (18:1.3) -- (54:1.7) -- (90:1.3);

\node at (90:0.45) {\small $\alpha$};
\node at (162:0.45) {\small $\beta$};
\node at (18:0.45) {\small $\gamma$};
\node at (234:0.45) {\small $\delta$};
\node at (-54:0.45) {\small $\epsilon$};

\node[draw,shape=circle, inner sep=1] at (0,0) {\small 1};

\foreach \x in {2,...,6}
\node[draw,shape=circle, inner sep=1] at (-90+72*\x:1) {\small \x};

\end{scope}
}


\draw[line width=1.5]
	(-18:1.7) -- (18:1.3) -- (18:0.7) -- (90:0.7) -- (162:0.7) -- (162:1.3) -- (198:1.7)
	(234:1.3) -- (-90:1.7) -- (-54:1.3);

\node at (30:0.8) {\small $\alpha$}; 
\node at (148:0.8) {\small $\alpha$}; 
\node at (174:0.8) {\small $\gamma$}; 
\node at (220:0.75) {\small $\epsilon$};
\node at (224:1.15) {\small $\delta$};
\node at (196:1.45) {\small $\beta$};
\node at (170:1.15) {\small $\alpha$};
\node at (-90:1.45) {\small $\alpha$}; 
\node at (2:0.75) {\small $\beta$}; 
\node at (-40:0.75) {\small $\delta$};
\node at (-46:1.15) {\small $\epsilon$};
\node at (-18:1.45) {\small $\gamma$};
\node at (10:1.15) {\small $\alpha$};

\node at (-54:1.8) {1};


\begin{scope}[xshift=3.5cm]

\draw[line width=1.5]
	(162:1.3) -- (198:1.7) -- (234:1.3) -- (-90:1.7) -- (-54:1.3) -- (-18:1.7) -- (18:1.3)
	(18:0.7) -- (90:0.7) -- (162:0.7)
	(90:0.7) -- (90:1.3);

\node at (76:0.8) {\small $\alpha$}; 
\node at (30:0.85) {\small $\beta$};
\node at (28:1.15) {\small $\delta$};
\node at (82:1.15) {\small $\gamma$};
\node at (54:1.5) {\small $\epsilon$};
\node at (104:0.8) {\small $\alpha$}; 
\node at (148:0.8) {\small $\gamma$};
\node at (154:1.2) {\small $\epsilon$};
\node at (99:1.1) {\small $\beta$};
\node at (128:1.5) {\small $\delta$};
\node at (174:0.8) {\small $\epsilon$}; 
\node at (220:0.75) {\small $\delta$};
\node at (226:1.1) {\small $\beta$};
\node at (198:1.45) {\small $\alpha$};
\node at (170:1.15) {\small $\gamma$}; 
\node at (-68:0.8) {\small $\delta$}; 
\node at (-116:0.8) {\small $\epsilon$};
\node at (-116:1.18) {\small $\gamma$};
\node at (-90:1.45) {\small $\alpha$};
\node at (-64:1.2) {\small $\beta$};
\node at (6:0.75) {\small $\epsilon$}; 
\node at (-40:0.75) {\small $\delta$};
\node at (-44:1.15) {\small $\beta$};
\node at (-18:1.45) {\small $\alpha$};
\node at (10:1.15) {\small $\gamma$};

\node at (-54:1.8) {2};

\end{scope}


\begin{scope}[xshift=7cm]

\draw[line width=1.5]
	(198:1.7) -- (234:1.3) -- (-90:1.7) 
	(234:1.3) -- (234:0.7)
	(-54:1.3) -- (-18:1.7) -- (18:1.3)
	(18:0.7) -- (90:0.7) -- (162:0.7)
	(90:0.7) -- (90:1.3);

\node at (76:0.8) {\small $\alpha$}; 
\node at (104:0.8) {\small $\alpha$}; 
\node at (226:1.15) {\small $\alpha$}; 
\node at (-116:1.15) {\small $\alpha$}; 
\node at (-18:1.45) {\small $\alpha$}; 

\node at (-54:1.8) {3};
	
\end{scope}

\end{tikzpicture}
\caption{Neighborhood tilings for $a^3b^2$.}
\label{3a2b}
\end{figure}

Next we turn to the edge length combination $a^3b^2$. By Proposition \ref{edge_combo} (also see \cite[Proposition 7]{gsy}), the edge lengths of the pentagon are arranged as the fourth of Figure \ref{edges1}. Up to the combinatorial symmetry of neighborhood, we may assume that the edges of $P_1$ are given by Figure \ref{3a2b}. We ignore the angles for the moment, and concentrate on the edge lengths. If the edge shared by $P_2,P_3$ is $a$, as in the first of Figure \ref{3a2b}, then we may successively determine all edges of $P_2,P_3,P_4,P_6,P_5$. If the edge shared by $P_2,P_3$ is $b$, as in the second and third of Figure \ref{3a2b}, then we may determine all edges of $P_2,P_3$. Since the two $b$-edges of $P_5$ are adjacent, either the edge shared by $P_4,P_5$ is $a$, or the edge shared by $P_5,P_6$ is $a$. Up to the symmetry of horizontal flipping, we may assume that $P_5,P_6$ share $a$. Then depending whether $P_4,P_5$ share $a$ or $b$, we get the second and third of Figure \ref{3a2b}. We label the three neighborhood tilings in Figure \ref{3a2b} (with angles to be discussed) as Cases 1, 2, 3.

Now we consider the angles. We denote the angles by $\alpha,\beta,\gamma,\delta,\epsilon$ as in $P_1$ in Figure \ref{3a2b}. First we argue that $\beta\ne\gamma$ and $\delta\ne\epsilon$. By Lemma \ref{geometry1} (or \cite[Lemma 21]{gsy}), we know that $\beta\ne\gamma$ is equivalent to $\delta\ne\epsilon$. Therefore we only need to assume $\beta=\gamma$ and $\delta=\epsilon$. In each of three pictures in Figure \ref{3a2b}, we have a vertex with three $a$ and another vertex with one $a$ and two $b$. The angle sums at the two vertices are $\alpha+2\beta=3\delta=2\pi$. This implies $\alpha+\beta+\gamma+\delta+\epsilon=\alpha+2\beta+2\delta=\frac{10}{3}\pi$. By the angle sum equation \eqref{asump} for pentagon, we conclude $f=12$. Therefore the case can be dismissed.

The case of $a^3b^2$ differs from $a^2b^2c$ and $a^3bc$ in that the angles cannot be uniquely described by the two bounding edges (with the exception of the $b^2$-angle $\alpha$). After showing $\beta\ne\gamma$ and $\delta\ne\epsilon$, the situation becomes less ambiguous.

To facilitate further discussion, we adopt the notation of \cite{gsy}.  Denote by $V_{ijk}$ the vertex shared by $P_i,P_j,P_k$. Denote by $A_{i,jk}$ the angle of $P_i$ at $V_{ijk}$. 

\subsubsection*{Case 1}

In the first of Figure \ref{3a2b}, we already know all edges and angles of $P_1$. We also know the $b^2$-angles $\alpha$ of $P_2,P_3$. Moreover, the $ab$-angles $A_{4,13},A_{6,12}$ are $\beta$ or $\gamma$. By the angle sums at $V_{134},V_{126}$ and $\beta\ne\gamma$, we know $A_{4,13}=\gamma$, $A_{6,12}=\beta$. This determines (all edges and angles of) $P_4,P_6$. 

Given we know all edges of $P_5$, there are two possible ways of arranging angles of $P_5$. Either way implies that one of $V_{145}$ and $V_{156}$ is $\delta^2\epsilon$, and the other is $\delta\epsilon^2$. The angle sums at $\delta^2\epsilon,\delta\epsilon^2$ and $V_{126}=\alpha\beta\gamma$ imply $\delta=\epsilon=\frac{2}{3}\pi$ and $\alpha+\beta+\gamma=2\pi$. Then we get $\alpha+\beta+\gamma+\delta+\epsilon=\frac{10}{3}\pi$. By the angle sum equation \eqref{asump} for pentagon, we conclude $f=12$. Therefore the case can be dismissed.

\subsubsection*{Case 2}

In the second of Figure \ref{3a2b}, we already know all edges and angles of $P_1$. The vertices $V_{145},V_{156}$ have only $a$-edges and therefore involve only $\delta,\epsilon$. By the angle sum at the two vertices and the fact that $\delta\ne\epsilon$, either $V_{145}=V_{156}=\delta^2\epsilon$, or $V_{145}=V_{156}=\delta\epsilon^2$. If $V_{145}=V_{156}=\delta^2\epsilon$, then by the fact that there is only one $\delta$ and one $\epsilon$ in $P_5$, the distribution of angles at the two vertices must be as indicated. This further determines all angles of $P_4,P_5,P_6$. If $V_{145}=V_{156}=\delta\epsilon^2$, we can similarly determine $P_4,P_5,P_6$. Since exchanging $\beta,\gamma$ and exchanging $\delta,\epsilon$ switches between the two scenarios, we only need to consider the case in the second of Figure \ref{3a2b}. The angle sums at $\alpha^3,\beta\gamma\epsilon,\delta^2\epsilon$ and the angle sum equation \eqref{asump} for pentagon imply
\begin{equation}\label{eq_case2}
\alpha=\tfrac{2}{3}\pi,\;
\beta+\gamma=\left(\tfrac{2}{3}+\tfrac{8}{f}\right)\pi,\;
\delta=\left(\tfrac{1}{3}+\tfrac{4}{f}\right)\pi,\;
\epsilon=\left(\tfrac{4}{3}-\tfrac{8}{f}\right)\pi.
\end{equation}

\subsubsection*{Case 3}

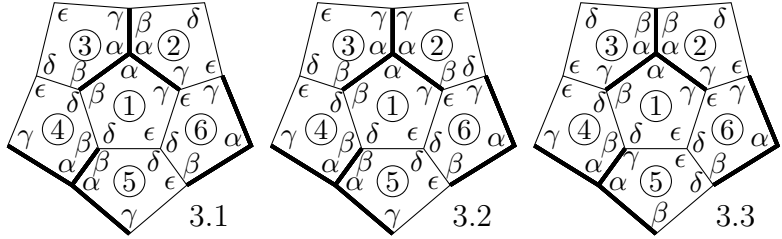
\begin{figure}[htp]
\centering
\begin{tikzpicture}[>=latex,scale=1]

\foreach \a in {0,1,2}
{
\begin{scope}[xshift=3.5*\a cm]

\foreach \x in {1,...,5}
{
\draw[rotate=72*\x]
	(90:0.7) -- (18:0.7) -- (18:1.3) -- (54:1.7) -- (90:1.3);
\draw[line width=1.5]
	(198:1.7) -- (234:1.3) -- (-90:1.7) 
	(234:1.3) -- (234:0.7)
	(-54:1.3) -- (-18:1.7) -- (18:1.3)
	(18:0.7) -- (90:0.7) -- (162:0.7)
	(90:0.7) -- (90:1.3);
}

\node at (90:0.45) {\small $\alpha$}; 
\node at (162:0.45) {\small $\beta$};
\node at (18:0.45) {\small $\gamma$};
\node at (234:0.45) {\small $\delta$};
\node at (-54:0.45) {\small $\epsilon$};
\node at (176:0.75) {\small $\delta$}; 
\node at (220:0.8) {\small $\beta$};
\node at (224:1.15) {\small $\alpha$};
\node at (198:1.45) {\small $\gamma$};
\node at (170:1.2) {\small $\epsilon$}; 

\node at (76:0.8) {\small $\alpha$}; 
\node at (104:0.8) {\small $\alpha$}; 
\node at (-116:1.15) {\small $\alpha$}; 
\node at (-18:1.45) {\small $\alpha$}; 

\node[draw,shape=circle, inner sep=1] at (0,0) {\small 1};

\foreach \x in {2,...,6}
\node[draw,shape=circle, inner sep=1] at (-90+72*\x:1) {\small \x};

\end{scope}
}


\node at (30:0.8) {\small $\gamma$};  
\node at (26:1.2) {\small $\epsilon$};
\node at (80:1.1) {\small $\beta$};
\node at (54:1.45) {\small $\delta$};

\node at (148:0.8) {\small $\beta$};  
\node at (152:1.2) {\small $\delta$};
\node at (99:1.15) {\small $\gamma$};
\node at (126:1.5) {\small $\epsilon$};

\node at (-66:0.8) {\small $\delta$}; 
\node at (-116:0.85) {\small $\beta$};
\node at (-62:1.15) {\small $\epsilon$};
\node at (-90:1.5) {\small $\gamma$};

\node at (6:0.75) {\small $\epsilon$}; 
\node at (-40:0.75) {\small $\delta$};
\node at (-44:1.15) {\small $\beta$};
\node at (8:1.1) {\small $\gamma$};

\node at (-54:1.8) {3.1};


\begin{scope}[xshift=3.5cm]

\node at (30:0.85) {\small $\beta$};  
\node at (28:1.15) {\small $\delta$};
\node at (82:1.15) {\small $\gamma$};
\node at (54:1.5) {\small $\epsilon$};

\node at (148:0.8) {\small $\beta$};  
\node at (152:1.2) {\small $\delta$};
\node at (99:1.15) {\small $\gamma$};
\node at (126:1.5) {\small $\epsilon$};

\node at (-66:0.8) {\small $\delta$}; 
\node at (-116:0.85) {\small $\beta$};
\node at (-62:1.15) {\small $\epsilon$};
\node at (-90:1.5) {\small $\gamma$};

\node at (6:0.75) {\small $\epsilon$}; 
\node at (-40:0.75) {\small $\delta$};
\node at (-44:1.15) {\small $\beta$};
\node at (8:1.1) {\small $\gamma$};

\node at (-54:1.8) {3.2};

\end{scope}


\begin{scope}[xshift=7cm]

\node at (30:0.8) {\small $\gamma$};  
\node at (26:1.2) {\small $\epsilon$};
\node at (80:1.1) {\small $\beta$};
\node at (54:1.45) {\small $\delta$};

\node at (148:0.8) {\small $\gamma$};  
\node at (154:1.2) {\small $\epsilon$};
\node at (99:1.1) {\small $\beta$};
\node at (126:1.45) {\small $\delta$};

\node at (-66:0.8) {\small $\epsilon$}; 
\node at (-114:0.8) {\small $\gamma$};
\node at (-62:1.15) {\small $\delta$};
\node at (-88:1.45) {\small $\beta$};

\node at (6:0.75) {\small $\epsilon$}; 
\node at (-40:0.75) {\small $\delta$};
\node at (-44:1.15) {\small $\beta$};
\node at (8:1.1) {\small $\gamma$};

\node at (-54:1.8) {3.3};

\end{scope}

\end{tikzpicture}
\caption{Neighborhood tilings for Case 3.}
\label{3a2bCase3}
\end{figure}

\begin{figure}[htp]
\centering
\begin{tikzpicture}[>=latex,scale=1]

\foreach \a in {0,1,2}
{
\begin{scope}[xshift=3.5*\a cm]

\foreach \x in {1,...,5}
{
\draw[rotate=72*\x]
	(90:0.7) -- (18:0.7) -- (18:1.3) -- (54:1.7) -- (90:1.3);
\draw[line width=1.5]
	(198:1.7) -- (234:1.3) -- (-90:1.7) 
	(234:1.3) -- (234:0.7)
	(-54:1.3) -- (-18:1.7) -- (18:1.3)
	(18:0.7) -- (90:0.7) -- (162:0.7)
	(90:0.7) -- (90:1.3);
}

\node at (90:0.45) {\small $\alpha$}; 
\node at (162:0.45) {\small $\beta$};
\node at (18:0.45) {\small $\gamma$};
\node at (234:0.45) {\small $\delta$};
\node at (-54:0.45) {\small $\epsilon$};
\node at (174:0.75) {\small $\epsilon$}; 
\node at (220:0.8) {\small $\gamma$};
\node at (224:1.15) {\small $\alpha$};
\node at (196:1.45) {\small $\beta$};
\node at (170:1.2) {\small $\delta$}; 

\node at (76:0.8) {\small $\alpha$}; 
\node at (104:0.8) {\small $\alpha$}; 
\node at (-116:1.15) {\small $\alpha$}; 
\node at (-18:1.45) {\small $\alpha$}; 

\node[draw,shape=circle, inner sep=1] at (0,0) {\small 1};

\foreach \x in {2,...,6}
\node[draw,shape=circle, inner sep=1] at (-90+72*\x:1) {\small \x};

\end{scope}
}


\node at (30:0.85) {\small $\beta$};  
\node at (28:1.2) {\small $\delta$};
\node at (80:1.15) {\small $\gamma$};
\node at (54:1.5) {\small $\epsilon$};

\node at (148:0.8) {\small $\beta$};  
\node at (152:1.2) {\small $\delta$};
\node at (99:1.15) {\small $\gamma$};
\node at (126:1.5) {\small $\epsilon$};

\node at (-66:0.8) {\small $\delta$}; 
\node at (-116:0.85) {\small $\beta$};
\node at (-62:1.15) {\small $\epsilon$};
\node at (-90:1.5) {\small $\gamma$};

\node at (6:0.75) {\small $\delta$}; 
\node at (-40:0.75) {\small $\epsilon$};
\node at (-44:1.15) {\small $\gamma$};
\node at (8:1.1) {\small $\beta$};


\begin{scope}[xshift=3.5cm]

\node at (30:0.8) {\small $\gamma$};  
\node at (26:1.2) {\small $\epsilon$};
\node at (80:1.1) {\small $\beta$};
\node at (54:1.45) {\small $\delta$};

\node at (-66:0.8) {\small $\epsilon$}; 
\node at (-114:0.8) {\small $\gamma$};
\node at (-62:1.15) {\small $\delta$};
\node at (-88:1.45) {\small $\beta$};

\node at (6:0.75) {\small $\delta$}; 
\node at (-40:0.75) {\small $\epsilon$};
\node at (-44:1.15) {\small $\gamma$};
\node at (8:1.1) {\small $\beta$};

\end{scope}


\begin{scope}[xshift=7cm]

\node at (30:0.85) {\small $\beta$};  
\node at (28:1.15) {\small $\delta$};
\node at (82:1.15) {\small $\gamma$};
\node at (54:1.5) {\small $\epsilon$};

\node at (148:0.8) {\small $\gamma$};  
\node at (154:1.2) {\small $\epsilon$};
\node at (99:1.1) {\small $\beta$};
\node at (128:1.5) {\small $\delta$};

\node at (-66:0.8) {\small $\epsilon$}; 
\node at (-114:0.8) {\small $\gamma$};
\node at (-62:1.15) {\small $\delta$};
\node at (-88:1.45) {\small $\beta$};

\node at (6:0.75) {\small $\epsilon$}; 
\node at (-40:0.75) {\small $\delta$};
\node at (-44:1.15) {\small $\beta$};
\node at (8:1.1) {\small $\gamma$};

\end{scope}

\end{tikzpicture}
\caption{More neighborhood tilings for Case 3.}
\label{3a2bCase3more}
\end{figure}
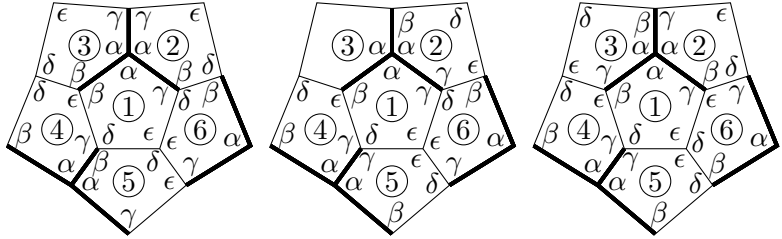

In the third of Figure \ref{3a2b}, we already know all edges. We also know all angles of $P_1$ and the $\alpha$ angles of all tiles. The angles of $P_4$ can be arranged in two ways. Figure \ref{3a2bCase3} gives the first way. Figure \ref{3a2bCase3more} gives the second way. 

In Figure \ref{3a2bCase3}, we know the $ab$-angles $A_{3,14},A_{5,14}$ are $\beta$ or $\gamma$. Moreover, the angle sums at $V_{134},V_{145}$ imply that $A_{3,14},A_{5,14}$ have the same value. The first and second of Figure \ref{3a2bCase3} are the case $A_{3,14}=A_{5,14}=\beta$, and the third picture is the case $A_{3,14}=A_{5,14}=\gamma$. The angles $A_{3,14},A_{5,14}$ determine $P_3,P_5$. 

The tile $P_6$ has two possible arrangements. All three in Figure \ref{3a2bCase3} have the same arrangement, and are labeled as Cases 3.1, 3.2, 3.3 (for $a^3b^2$). Then in the first two pictures, we further consider the two possible arrangements of $P_2$. In the first picture, the angle sums at $\alpha^3,\beta^2\delta,\gamma^2\epsilon,\delta^2\epsilon$ and the angle sum equation \eqref{asump} for pentagon imply
\begin{equation}\label{eq_case31}
\alpha=\tfrac{2}{3}\pi,\;
\beta=\left(\tfrac{1}{3}+\tfrac{4}{f}\right)\pi,\;
\gamma=\delta=\left(\tfrac{4}{3}-\tfrac{8}{f}\right)\pi,\;
\epsilon=\left(-\tfrac{2}{3}+\tfrac{16}{f}\right)\pi.
\end{equation}
In the second picture, the angle sums at $\alpha^3,\beta^2\delta,\beta\gamma\epsilon,\delta^2\epsilon$ and the angle sum equation \eqref{asump} for pentagon imply
\begin{equation}\label{eq_case32}
\alpha=\tfrac{2}{3}\pi,\;
\beta=\left(\tfrac{5}{6}-\tfrac{2}{f}\right)\pi,\;
\gamma=\left(-\tfrac{1}{6}+\tfrac{10}{f}\right)\pi,\;
\delta=\left(\tfrac{1}{3}+\tfrac{4}{f}\right)\pi,\;
\epsilon=\left(\tfrac{4}{3}-\tfrac{8}{f}\right)\pi.
\end{equation}
In the third picture, we also consider two possible arrangements of $P_2$. If the arrangement is not the same as the third picture, then the angle sums at $V_{134}=\beta\gamma\delta$ and $V_{126}=\beta\gamma\epsilon$ imply that $\delta=\epsilon$, a contradiction. Therefore the angles of $P_2$ must be arranged as indicated. Then the angle sums at $\alpha^3,\beta\gamma\delta,\gamma^2\epsilon,\delta\epsilon^2$ and the angle sum equation \eqref{asump} for pentagon imply
\begin{equation}\label{eq_case33}
\alpha=\tfrac{2}{3}\pi,\;
\beta=\left(-\tfrac{1}{6}+\tfrac{10}{f}\right)\pi,\;
\gamma=\left(\tfrac{5}{6}-\tfrac{2}{f}\right)\pi,\;
\delta=\left(\tfrac{4}{3}-\tfrac{8}{f}\right)\pi,\;
\epsilon=\left(\tfrac{1}{3}+\tfrac{4}{f}\right)\pi.
\end{equation}
We remark that by $f\ge 16$, we have $\beta<\gamma$ and $\delta>\epsilon$ for \eqref{eq_case31} (Case 3.1) and \eqref{eq_case33} (Case 3.3), and we also have $\beta>\gamma$ and $\delta<\epsilon$ for \eqref{eq_case32} (Case 3.2). Therefore the geometric constraint in Lemma \ref{geometry1} is satisfied for all three cases.

We still need to consider the other arrangement of $P_6$, given we already know $P_1,P_3,P_4,P_5$. The new arrangement of $P_6$ means that $V_{126}=\beta\gamma\delta$ or $\gamma^2\delta$, and $V_{156}=\epsilon^3$. In the first two pictures of Figure \ref{3a2bCase3}, we compare the angle sums at $V_{134}=\beta^2\delta$ and $V_{126}$, and always get $\beta=\gamma$, a contradiction. In the third picture, the angle sums at $\alpha^3,\beta\gamma\delta,\epsilon^3$ imply $\alpha=\epsilon=\frac{2}{3}\pi$ and $\beta+\gamma+\delta=2\pi$. Then we get $\alpha+\beta+\gamma+\delta+\epsilon=\frac{10}{3}\pi$. By the angle sum equation \eqref{asump} for pentagon, we conclude $f=12$, and the case can be dismissed.

Now we turn to Figure \ref{3a2bCase3more}, where the angles of $P_4$ are arranged differently from Figure \ref{3a2bCase3}. Then we consider two possible arrangements of $P_5$. The first picture shows one arrangement, and the second and third pictures show the other arrangement. 

In the first of Figure \ref{3a2bCase3more}, we have $A_{3,14}=\beta$ or $\gamma$. If $A_{3,14}=\gamma$, then the angle sums at $V_{134}=\beta\gamma\epsilon$ and $V_{145}=\beta\gamma\delta$ imply $\delta=\epsilon$, a contradiction. Therefore $A_{3,14}=\beta$, which determines $P_3$. Now we consider two possible arrangements of $P_6$. If $P_6$ is arranged differently from the first of Figure \ref{3a2bCase3more}, then $V_{126}=\beta\gamma\epsilon$ or $\gamma^2\epsilon$. Comparing the angle sums at $V_{134}=\beta^2\epsilon$ and $V_{126}$ always gives $\beta=\gamma$, a contradiction. This implies that $P_6$ must be arranged as in the picture. Moreover, comparing the angle sums at $V_{145}=\beta\gamma\delta$ and $V_{126}$ gives $A_{2,16}=\beta$. This determines $P_2$. Now the neighborhood tiling has a collection of vertices $\alpha^3,\beta^2\epsilon,\beta\gamma\delta,\delta\epsilon^2$. The collection is obtained from the collection for Case 2 by exchanging $\delta$ and $\epsilon$. Then the equalities $\beta>\gamma$ and $\delta<\epsilon$ for Case 2 become $\beta>\gamma$ and $\delta>\epsilon$ for the first of Figure \ref{3a2bCase3more}. This contradicts the geometric constraint in Lemma \ref{geometry1}.

In the second and third pictures of Figure \ref{3a2bCase3more}, we already know $P_1,P_4,P_5$ and then consider two possible arrangements of $P_6$. The second picture shows the first arrangement of $P_6$. By comparing the angle sums at $V_{145}=\gamma^2\delta$ and $V_{126}$, we find $A_{2,16}=\gamma$, which determines $P_2$. Then we have two possible arrangements of $P_3$. One arrangement of $P_3$ has vertices $\alpha^3,\beta^2\epsilon,\gamma^2\delta,\epsilon^3$, and the other arrangement has vertices $\alpha^3,\beta\gamma\epsilon,\gamma^2\delta,\epsilon^3$. Combined with the angle sum equation \eqref{asump} for pentagon, we get
\[
\alpha=\beta=\epsilon=\tfrac{2}{3}\pi,\;
\gamma=\left(1-\tfrac{4}{f}\right)\pi,\;
\delta=\tfrac{8}{f}\pi,
\]
and
\[
\alpha=\epsilon=\tfrac{2}{3}\pi,\;
\beta=\left(\tfrac{1}{2}+\tfrac{2}{f}\right)\pi,\;
\gamma=\left(\tfrac{5}{6}-\tfrac{2}{f}\right)\pi,\;
\delta=\left(\tfrac{1}{3}+\tfrac{4}{f}\right)\pi.
\]
By $f\ge 16$, we have $\beta<\gamma$ and $\delta<\epsilon$ in both cases, contradicting the geometric constraint in Lemma \ref{geometry1}. 

In the third of Figure \ref{3a2bCase3more}, we have the other arrangement of $P_6$. By comparing the angle sums at $V_{126}$ and $V_{134}$, we find $A_{2,16}=\beta,A_{3,14}=\gamma$. This determines $P_2,P_3$. The neighborhood tiling has vertices $\alpha^3,\beta\gamma\epsilon,\gamma^2\delta,\delta\epsilon^2$. The angle sums at the vertices imply $\beta=\delta$ and $\gamma=\epsilon$, contradicting the geometric constraint in Lemma \ref{geometry1}. 

The following summarises the discussion on tilings of the neighborhood of a $3^5$-tile with edge combination $a^3b^2$.

\begin{proposition}\label{main3}
If an edge-to-edge spherical tiling by more than $12$ congruent pentagons has edge length combination $a^3b^2$, with $a,b$ distinct, then up to symmetry, the neighborhood of a $3^5$-tile has four possible tilings given by the second of Figure \ref{3a2b} and the three pictures in Figure \ref{3a2bCase3}. 
\end{proposition}

\section{Tiling for Edge Combination $a^3b^2$}
\label{3a2b_tiling}

We need to start from the four neighborhood tilings in Proposition \ref{main3} and try to extend the tiling beyond the neighborhood. The following technical results are useful for constructing the extension.

\begin{lemma}\label{beven}
In a tiling by congruent pentagons with edge length combination $a^3b^2$, $a,b$ distinct, the number of $ab$-angles at any vertex is even. 
\end{lemma}

\begin{proof}
The proposition is a consequence of the fact that the number of $b$-edges at a vertex is the number of $b^2$-angles plus half of the number of $ab$-angles.
\end{proof}

\begin{lemma}\label{klem}
Suppose $\beta\ne\gamma$ and $\delta\ne\epsilon$ in the $a^3b^2$-pentagon in Figure \ref{geom7}. Suppose a tiling by congruent copies of the pentagon has at most one $\epsilon$ at every vertex. Then there cannot be consecutive $\gamma\delta\cdots\delta\gamma$ at a vertex ($\cdots$ consists of $\delta$ only). In case there is no $\delta$, this means that two $\gamma$ cannot share an $a$-edge at a vertex.
\end{lemma}

We note that the lemma is symmetric with respect to the exchange of $\beta,\gamma$ and exchange of $\delta,\epsilon$.

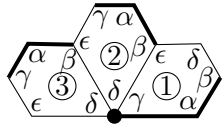
\begin{figure}[htp]
\centering
\begin{tikzpicture}[>=latex,scale=1]

\draw
	(0,0) -- (60:1.1) -- (40:1.5) -- (20:1.5)
	(100:1.5) -- (120:1.1) -- (0,0)
	(160:1.5) -- (180:1.1) -- (0,0);

\draw[line width=1.5]
	(0,0) -- (0:1.1) -- (20:1.5)
	(60:1.1) -- (80:1.5) -- (100:1.5)
	(120:1.1) -- (140:1.5) -- (160:1.5);

\fill (0,0) circle (0.1);

\node at (30:0.35) {\small $\gamma$}; 
\node at (22:1.3) {\small $\beta$};
\node at (10:1) {\small $\alpha$};
\node at (36:1.3) {\small $\delta$};
\node at (52:1) {\small $\epsilon$}; 

\node at (90:0.35) {\small $\delta$}; 
\node at (112:1.05) {\small $\epsilon$};
\node at (84:1.3) {\small $\alpha$};
\node at (98:1.3) {\small $\gamma$};
\node at (70:0.95) {\small $\beta$};

\node at (150:0.35) {\small $\delta$}; 
\node at (172:1.05) {\small $\epsilon$};
\node at (142:1.3) {\small $\alpha$};
\node at (158:1.3) {\small $\gamma$};
\node at (130:0.95) {\small $\beta$};

\node[draw,shape=circle, inner sep=0.5] at (30:0.8) {\small $1$};
\node[draw,shape=circle, inner sep=0.5] at (90:0.8) {\small $2$};
\node[draw,shape=circle, inner sep=0.5] at (150:0.8) {\small $3$};

\end{tikzpicture}
\caption{Cannot have consecutive $\gamma\delta^k\gamma$ if $\epsilon^2\cdots$ is not a vertex.}
\label{impossible}
\end{figure}

\begin{proof}
Suppose we have consecutive $\gamma\delta\cdots\delta$ at a vertex as in Figure \ref{impossible}. Then we may determine all edges and angles of $P_1$. Since $\delta$ is adjacent to $\beta$ and $\epsilon$ in the pentagon, $P_1,P_2$ share a vertex $\beta\epsilon\cdots$ or $\epsilon^2\cdots$. Since $\epsilon^2\cdots$ is not a vertex, the vertex shared by $P_1,P_2$ is $\beta\epsilon\cdots$. This determines $\beta$ of $P_2$. Combined with $\delta$ of $P_2$, we determine $P_2$. By applying the similar argument to the vertex shared by $P_2,P_3$, we also determine $P_3$. The process continues.

If we have consecutive $\gamma\delta\cdots\delta\gamma$, then we may start from one end $\gamma$ and consecutively determine the tiles containing the intermediate $\delta$. When we reach the other end $\gamma$, we find that the two tiles containing consecutive $\gamma$ and $\delta$ share a vertex $\epsilon^2\cdots$, contradicting to the assumption. 
\end{proof}

\subsection{Case 2, $\alpha=\beta$}
\label{type2A}

The neighborhood tiling for Case 2 is the second of Figure \ref{3a2b}, and the angles are given by \eqref{eq_case2}. By $f\ge 16$, we have $\epsilon>\alpha>\delta$ and $\beta+\gamma<2\alpha$. By Lemma \ref{geometry1}, we also have $\beta>\gamma$. Then by $\beta+\gamma<2\alpha$ and $\beta+\gamma=2\delta$, we get $\gamma<\alpha$ and $\delta\ne \beta,\gamma$. By $\beta+\gamma<\epsilon$, we also have $\epsilon\ne \beta,\gamma$. We conclude that the only way for some of the five angles to be equal is $\alpha=\beta$. This section studies the case $\alpha=\beta$, and the next section studies the case $\alpha\ne\beta$. 

By $\alpha=\beta$ and \eqref{eq_case2}, we get all the angles
\[
\alpha=\beta=\frac{2}{3}\pi, \;
\gamma=\frac{8}{f}\pi, \;
\delta=\left(\frac{1}{3}+\frac{4}{f}\right)\pi, \;
\epsilon=\left(\frac{4}{3}-\frac{8}{f}\right)\pi.
\]
By $f\ge 16$, we have $2\pi-2\epsilon=\left(-\frac{2}{3}+\frac{16}{f}\right)\pi<$ all angles. This implies that $\epsilon^2\cdots$ is not a vertex, and therefore Lemma \ref{klem} can be applied.

In the second of Figure \ref{3a2b}, $P_5,P_6$ share a vertex $\beta^2\cdots$, with angles in the remainder $\cdots$ being consecutive. Since the angle sum of the remainder is $2\pi-2\beta=\frac{2}{3}\pi=\alpha=\beta<2\delta,\epsilon$, by Lemma \ref{beven}, the vertex $\beta^2\cdots$ is $\alpha\beta^2$, $\beta^2\gamma^c\delta$, or $\beta^2\gamma^c$. By Lemma \ref{klem}, the remainder cannot be $\gamma^c\delta$ or $\gamma^c$. Therefore the vertex $\beta^2\cdots$ shared by $P_5,P_6$ is $\alpha\beta^2$. This $\alpha$ belongs to a tile outside $P_5,P_6$. Since $\gamma$ is adjacent to $\alpha$ in the pentagon, the new tile shares a vertex $\alpha\gamma\cdots$ with either $P_5$ or $P_6$. 

By Lemma \ref{beven}, the vertex $\alpha\gamma\cdots$ is either $\alpha\beta\gamma\cdots$ or $\alpha\gamma^2\cdots$. If $\alpha\beta\gamma\cdots$ is a vertex, then the remainder has angle sum $\left(\frac{2}{3}-\frac{8}{f}\right)\pi<\alpha,\beta,2\delta,\epsilon$. Therefore the vertex is $\alpha\beta\gamma^c\delta$ or $\alpha\beta\gamma^c$, $c\ge 1$. If $c\ge 3$, then by the edge length consideration, we get two $\gamma$ sharing an $a$-edge, a contradiction to Lemma \ref{klem}. Therefore by Lemma \ref{beven}, the vertex is $\alpha\beta\gamma\delta$ or $\alpha\beta\gamma$. By $f>16$, we find the angle sum of $\alpha\beta\gamma$ to be $<2\pi$. Therefore $\alpha\beta\gamma\cdots$ can only be $\alpha\beta\gamma\delta$.

Now we turn to the possibility that $\alpha\gamma\cdots$ is $\alpha\gamma^2\cdots$. Having discussed $\alpha\beta\gamma\cdots$, we may assume that the remainder of $\alpha\gamma^2\cdots$ does not contain $\beta$. Since $\alpha+\gamma+\epsilon=2\pi$, the remainder does not contain $\epsilon$. Therefore the vertex is $\alpha\gamma^c\delta^d$ or $\alpha^2\gamma^c\delta^d$, with $c\ge 2$. By the edge length consideration, we get a contradiction to Lemma \ref{klem}.

We conclude that the vertex $\alpha\gamma\cdots$ shared by $P_5,P_6$ is $\alpha\beta\gamma\delta$. Then the angle sum at $\alpha\beta\gamma\delta$ implies $f=36$ and
\[
\alpha=\beta=\frac{2}{3}\pi, \;
\gamma=\frac{2}{9}\pi, \;
\delta=\frac{4}{9}\pi, \;
\epsilon=\frac{10}{9}\pi.
\]
Now consider the vertex $\beta\gamma\cdots$ shared by $P_4,P_5$ in the second of Figure \ref{3a2b}. The remainder has angle sum $\frac{10}{9}\pi$. By the angle length consideration and the angle sum estimation, the remainder cannot contain $\beta$ (which by Lemma \ref{beven} contains another $\beta$ or $\gamma$) or $\epsilon$. Therefore the vertex is $\beta\gamma^c\delta^d$ or $\alpha\beta\gamma^c\delta^d$. Since $\beta$ and $\gamma$ share $a$-edge, by the edge length consideration, we must have $c\ge 3$. However, this contradicts Lemma \ref{klem}.

\subsection{Case 2, $\alpha\ne\beta$}
\label{type2B}

We continue the study of Case 2, under the additional assumption that all five angles are distinct. From the discussion at the beginning of Section \ref{type2A}, we also know $\beta>\gamma$ and $\delta<\epsilon$.

Suppose $\epsilon^2\cdots$ is a vertex. Then $0>2\pi-2\epsilon=2\left(\frac{8}{f}-\frac{1}{3}\right)\pi$. This implies $f<24$. By \eqref{vcountf}, we get $5\ge v_4+2v_5+3v_6+\cdots$. By \cite[Theorems 1 and 6]{yan}, we get $v_{>6}=0$ for $f=22$, $v_{>5}=0$ for $f=20$, and $v_{>4}=0$ for $f=18$ or $16$. 

By $f\ge 16$, we have $2\pi-2\epsilon<\left(\frac{1}{3}+\frac{4}{f}\right)\pi\le \alpha,\beta,\delta,\epsilon$. Therefore the vertex $\epsilon^2\cdots$ must be $\gamma^k\epsilon^2$. By Lemma \ref{beven} and the fact that all vertices have degree $\le 6$, we get $\epsilon^2\cdots=\gamma^2\epsilon^2$ or $\gamma^4\epsilon^2$.

If $\gamma^4\epsilon^2$ is a vertex, then $v_6\ge 1$. By $5\ge v_4+2v_5+3v_6$ and \cite[Theorems 1 and 6]{yan}, we get $v_4=2,v_5=0,v_6=1$. By $v_{>6}=0$ and \eqref{vcountf}, we get $f=22$. By $f=22$ and the angle sum at $\gamma^4\epsilon^2$, we may calculate all the angles
\[
\alpha=\frac{2}{3}\pi,\;
\beta=\frac{67}{66}\pi,\;
\gamma=\frac{1}{66}\pi,\;
\delta=\frac{17}{33}\pi,\;
\epsilon=\frac{32}{33}\pi.
\]
Since no four angles from above (repetition allowed) can add up to $2\pi$, we get a contradiction to $v_4=2$. 

\begin{figure}[htp]
\centering
\begin{tikzpicture}[>=latex]

\foreach \a in {-1,0,1,2}
\draw[rotate=72*\a]
	(0,0) -- (90:0.8) -- (72:1.5) -- (36:1.5) -- (18:0.8) -- cycle;
	
\draw[line width=1.5]
	(0,0) -- (90:0.8)
	(0,0) -- (18:0.8) -- (0:1.5)
	(0,0) -- (162:0.8) -- (180:1.5);
	
\foreach \a in {1,...,4}
\node[draw,shape=circle,inner sep=1] at (-90+72*\a:0.95) {\small \a};

\foreach \a in {1,-1}
{
\begin{scope}[xscale=\a]

\node at (54:0.3) {\small $\alpha$};
\node at (75:0.75) {\small $\beta$};
\node at (33:0.75) {\small $\gamma$};
\node at (68:1.3) {\small $\delta$};
\node at (40:1.3) {\small $\epsilon$};

\node at (-20:0.3) {\small $\gamma$};
\node at (6:0.75) {\small $\alpha$};
\node at (-42:0.75) {\small $\epsilon$};
\node at (-8:1.3) {\small $\beta$};
\node at (-30:1.3) {\small $\delta$};

\end{scope}
}

\node at (-90:0.3) {\small $\delta$};

\end{tikzpicture}
\caption{Impossible vertex for Case 2 .}
\label{case2fig1}
\end{figure}
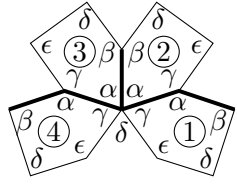

If $\gamma^2\epsilon^2$ is a vertex, then the angle sum at $\gamma^2\epsilon^2$ gives $\beta=\pi$ and $\gamma=\left(\frac{8}{f}-\frac{1}{3}\right)\pi$. For each $f=16,18,20,22$, we know the values of all five angles and can further find all possible angle combinations at vertices. We also take into account of Lemma \ref{beven} and the constraint $v_{>6}=0$ for $f=22$, $v_{>5}=0$ for $f=20$, and $v_{>4}=0$ for $f=18$ or $16$. Then we get the following {\em anglewise vertex combinations} (all possible angle combinations at vertices)
\begin{align*}
\text{AVC}
&=\{\alpha^3, \;
\beta\gamma\epsilon, \;
\delta^2\epsilon, \;
\gamma^2\epsilon^2\}
& \text{for }f
&=16,18,22; \\
\text{AVC}
&=\{\alpha^3, \;
\beta\gamma\epsilon, \;
\delta^2\epsilon, \;
\gamma^2\epsilon^2,  \;
\alpha^2\gamma^2\delta\}
&\text{for }f
&=20. 
\end{align*}
By the AVC, the vertex $\gamma\delta\cdots$ shared by $P_2,P_6$ in the second of Figure \ref{3a2b} is $\alpha^2\gamma^2\delta$, and we must have $f=20$. It is easy to see that $\alpha^2\gamma^2\delta$ must be configured as in Figure \ref{case2fig1}. This determines $P_1,P_4$. Then $P_1,P_2$ share a vertex $\alpha\beta\cdots$ or $\alpha\gamma\cdots$. Since $\alpha\beta\cdots$ is not in the AVC, we get the $\gamma$ angle of $P_2$, which further determines $P_2$. By the same argument, we determine $P_3$. Then $P_2,P_3$ share a vertex $\beta^2\cdots$, contradicting the AVC. This completes the proof that $\epsilon^2\cdots$ is not a vertex. Therefore Lemma \ref{klem} holds.

Next we find vertices $\alpha^a\beta^b\gamma^c\delta^d\epsilon$ with one $\epsilon$. Let $\bar{b}=\min\{b,c\}$. Then the angle sum at such a vertex gives
\begin{align}
2&\ge 
\frac{2}{3}a
+\left(\frac{2}{3}+\frac{8}{f}\right)\bar{b}
+\left(\frac{1}{3}+\frac{4}{f}\right)d
+\left(\frac{4}{3}-\frac{8}{f}\right) \nonumber \\
&=\frac{1}{3}(2a+2\bar{b}+d+4)+\frac{4}{f}(2\bar{b}+d-2). \label{case2eq2}
\end{align}
If $2\bar{b}+d-2>0$, then $2a+2\bar{b}+d<2$. This implies $a=\bar{b}=0$ and $d=0$ or $1$, and further implies $2\bar{b}+d-2\le 0$, contradicting to the assumption. Therefore $2\bar{b}+d\le 2$. 

If $2\bar{b}+d=2$, then \eqref{case2eq2} implies $a=0$. By $2\bar{b}+d=2$, we get vertices $\beta^b\gamma\epsilon$ ($b\ge 1$), $\beta\gamma^c\epsilon$ ($c\ge 1$), $\beta^b\delta^2\epsilon$, $\gamma^c\delta^2\epsilon$. In the second of Figure \ref{3a2b}, we see that $\beta\gamma\epsilon$ and $\delta^2\epsilon$ are already vertices. Then the angle sum implies that the four possibilities $\beta^b\gamma\epsilon$, $\beta\gamma^c\epsilon$, $\beta^b\delta^2\epsilon$, $\gamma^c\delta^2\epsilon$ can only be the existing vertices $\beta\gamma\epsilon$ and $\delta^2\epsilon$. 

If $2\bar{b}+d=1$, then \eqref{case2eq2} and $f\ge 16$ imply $a=0$. By $2\bar{b}+d=1$, we get vertices $\beta^b\delta\epsilon$ ($b\ge 1$), $\gamma^c\delta\epsilon$ ($c\ge 1$). By $\beta+\delta+\epsilon\ge \frac{1}{2}\left(\frac{2}{3}+\frac{8}{f}\right)\pi+\left(\frac{1}{3}+\frac{4}{f}\right)\pi+\left(\frac{4}{3}-\frac{8}{f}\right)\pi>2\pi$, we can only have $\gamma^c\delta\epsilon$. By Propositions \ref{beven} and \ref{klem}, we further find that the only possibility is $\gamma^2\delta\epsilon$. 

If $2\bar{b}+d=0$, then \eqref{case2eq2} and $f\ge 16$ imply $a=0$ or $1$. By $2\bar{b}+d=0$, we get vertices $\alpha\beta^b\epsilon$ ($b\ge 1$), $\alpha\gamma^c\epsilon$ ($c\ge 1$), $\beta^b\epsilon$ ($b\ge 2$), $\gamma^c\epsilon$ ($c\ge 2$). By $\alpha+\beta+\epsilon>2\pi$ and $2\beta+\epsilon>2\pi$ are $>2\pi$, we can only have $\alpha\gamma^c\epsilon,\gamma^c\epsilon$. By Propositions \ref{beven} and \ref{klem}, we further find that the only possibilities are $\alpha\gamma^2\epsilon,\gamma^2\epsilon$. The angle sum at $\gamma^2\epsilon$ implies $\beta=\gamma$, a contradiction. The angle sum at $\alpha\gamma^2\epsilon$ implies $\beta=\left(\frac{2}{3}+\frac{4}{f}\right)\pi$ and $\gamma=\frac{4}{f}\pi$. Note that the vertex $\beta^2\cdots$ shared by $P_5,P_6$ in the second of Figure \ref{3a2b} has the remaining angle $2\pi-2\beta=\left(\frac{2}{3}-\frac{8}{f}\right)\pi$, which is strictly less than $\alpha,\beta,\epsilon$. Therefore the vertex is $\beta^2\gamma^c\delta^d$. The second of Figure \ref{3a2b} further shows that $\gamma^c\delta^d$ is bounded by two $b$-edges, contradicting Lemma \ref{klem}. 

It remains to consider a vertex $\alpha^a\beta^b\gamma^c\delta^d$ without $\epsilon$. To avoid contradicting Lemma \ref{klem}, each $\gamma$ needs to be combined with one $\beta$ into a chain $\beta\delta\cdots\delta\gamma$ bordered by two $b$-edges. This implies $b\ge c$. All the possible vertices we found so far give the following complete list of possible angle combinations at vertices
\begin{equation}\label{case2eq3}
\text{AVC} 
=\{\alpha^3, \;
\beta\gamma\epsilon, \;
\delta^2\epsilon, \;
\gamma^2\delta\epsilon, \;
\alpha^a\beta^b\gamma^c\delta^d(b\ge c)\}.
\end{equation}

Suppose $\gamma^2\delta\epsilon$ is not a vertex. Since both numbers of $\beta$ and $\gamma$ in the whole tiling should be equal to $f$, to balance the equal total number, we must have $b=c$ in every vertex $\alpha^a\beta^b\gamma^c\delta^d$. The angle sum at $\alpha^a\beta^b\gamma^b\delta^d$ is
\[
2=\dfrac{2}{3}a+\left(\frac{2}{3}+\frac{8}{f}\right)b+\left(\frac{1}{3}+\frac{4}{f}\right)d
=\dfrac{2}{3}a+\left(\frac{1}{3}+\frac{4}{f}\right)(2b+d).
\]
If the vertex is not $\alpha^3$, this implies $2a+2b+d\le 5$. By trying all $a,b,d$ satisfying the inequality and using $f\ge 16$, we get all the solutions
\begin{align*}
a=0, \; 2b+d=4, \; f=24 &\colon 
\beta\gamma\delta^2, \;
\beta^2\gamma^2, \;
\delta^4; \\
a=1, \; 2b+d=3, \; f=36 &\colon 
\alpha\beta\gamma\delta, \;
\alpha\delta^3; \\
a=0, \; 2b+d=5, \; f=60 &\colon 
\beta^2\gamma^2\delta, \;
\beta\gamma\delta^3, \;
\delta^5.
\end{align*}
The vertex $\beta^2\cdots$ shared by $P_5,P_6$ in the second of Figure \ref{3a2b} does not appear in the list for $f=36$, is $\beta^2\gamma^2$ for $f=24$, and is $\beta^2\gamma^2\delta$ for $f=60$. The second of Figure \ref{3a2b} further shows that $\gamma^2$ (for $f=24$) or $\gamma^2\delta$ (for $f=60$) is bounded by two $b$-edges, contradicting Lemma \ref{klem}. 

We conclude that $\gamma^2\delta\epsilon$ must be a vertex. The angle sum at $\gamma^2\delta\epsilon$ implies $\beta=\left(\frac{1}{2}+\frac{6}{f}\right)\pi$ and $\gamma=\left(\frac{1}{6}+\frac{2}{f}\right)\pi$. Then the angle sum at $\alpha^a\beta^b\gamma^c\delta^d$ ($b\ge c$) is
\[
\frac{2}{3}a
+\left(\frac{1}{6}+\frac{2}{f}\right)(3b+c+2d)
=2.
\]
If $3b+c+2d=0$, then the vertex is $\alpha^3$. If $3b+c+2d>0$, then the equation implies $4a+3b+c+2d<12$. Substituting those $a,3b+c+2d$ satisfying the inequality and that $3b+c+2d$ is positive and even (by Lemma \ref{beven}) into the equation and solve for $f$. Those combinations $(a,3b+c+2d)$ yielding even $f\ge 16$ are $(0,8)$ for $f=24$, $(1,6)$ for $f=36$, and $(0,10)$ for $f=60$.

For $f=24$, we consider the vertex $\beta^2\cdots$ shared by $P_5,P_6$ in the second of Figure \ref{3a2b}. By the AVC \eqref{case2eq3} and the discussion above, the vertex is $\beta^b\gamma^c\delta^d$, satisfying $b\ge 2$, $b\ge c$, $3b+c+2d=8$, and $b,c$ have the same parity. The solution shows that the vertex $\beta^2\cdots$ is $\beta^2\gamma^2$ or $\beta^2\delta$. The second of Figure \ref{3a2b} further shows that $\gamma^2$ or $\delta$ is bounded by two $b$-edges. By Lemma \ref{klem}, we cannot have $\gamma^2$ bounded by two $b$-edges. Moreover, $\delta$ is not a $b^2$-angle.

For $f=36$, we consider the vertex $\beta\gamma\cdots$ shared by $P_4,P_5$ in the second of Figure \ref{3a2b}. By the AVC \eqref{case2eq3} and the discussion above, the vertex is $\beta\gamma\epsilon$ or $\alpha\beta^b\gamma^c\delta^d$, satisfying $b\ge 1$, $c\ge 1$, $b\ge c$, $3b+c+2d=6$, and $b,c$ have the same parity. The solution shows that the vertex $\beta\gamma\cdots$ is $\beta\gamma\epsilon$ or $\alpha\beta\gamma\delta$. The second of Figure \ref{3a2b} further shows that $\epsilon$ or $\alpha\delta$ is bounded by two $b$-edges. Both are contradictions.

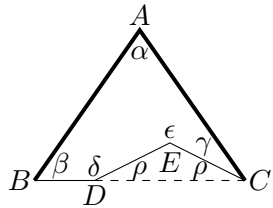
\begin{figure}[htp]
\centering
\begin{tikzpicture}

\draw[line width=1.5]
	(-1.4,0) -- (0,2) -- (1.4,0);
\draw
	(-1.4,0) -- (-0.6,0) -- (0.4,0.5) -- (1.4,0);
		
\draw[dashed]
	(-0.6,0) -- (1.4,0);

\node at (0,1.7) {\small $\alpha$};
\node at (-1.05,0.2) {\small $\beta$};
\node at (0.85,0.45) {\small $\gamma$};
\node at (-0.6,0.2) {\small $\delta$};
\node at (0.4,0.65) {\small $\epsilon$};

\node at (0,0.13) {\small $\rho$};
\node at (0.8,0.13) {\small $\rho$};
 
\node at (0,2.2) {\small $A$};
\node at (-1.6,0) {\small $B$};
\node at (1.6,0) {\small $C$};
\node at (-0.6,-0.2) {\small $D$};
\node at (0.4,0.25) {\small $E$};
           
\end{tikzpicture}
\caption{Impossible spherical pentagon.}
\label{case2fig2}
\end{figure}

For $f=60$, we have
\[
\alpha=\frac{2}{3}\pi,\;
\beta=\frac{3}{5}\pi,\;
\gamma=\frac{1}{5}\pi,\;
\delta=\frac{2}{5}\pi,\;
\epsilon=\frac{6}{5}\pi.
\]
Substituting the angles into \eqref{eq11}, we get 
\[
P=0,\quad
Q=2\left(\frac{\sqrt{10-2\sqrt{5}}}{4}\right)^3,\quad
R=\frac{5-2\sqrt{5}}{4}.
\]
Therefore
\[
\cot\theta=-\frac{R}{Q}=-\frac{\sqrt{5}-1}{\sqrt{10+2\sqrt{5}}}.
\]
Since $\alpha<\pi$, we have $\theta<\pi$. Therefore $\theta=\frac{3}{5}\pi=\beta$, and the pentagon is given by Figure \ref{case2fig2}, in which $D$ lies in the great arc connecting $B$ and $C$.

Since $\triangle CDE$ and $\triangle ABC$ are isosceles triangles, we have $\rho=\angle EDC=\angle ECD=\beta-\gamma=\frac{2}{5}\pi$. Then $\delta+\rho=\frac{4}{5}\pi\ne\pi$, contradicting to the fact that $D$ lies on the great arc connecting $B$ and $C$.

\subsection{Case 3.1}

The neighborhood tiling for Case 3.1 is the first of Figure \ref{3a2bCase3}, and the angles are given by \eqref{eq_case31}. For $\epsilon>0$, we must have $f<24$. This means we only need to consider $f=16,18,20,22$. 

For $f=16$, the angles are 
\[
\alpha=\frac{2}{3}\pi,\;
\beta=\frac{7}{12}\pi,\;
\gamma=\delta=\frac{5}{6}\pi,\;
\epsilon=\frac{1}{3}\pi.
\]
Consider the vertex  $\beta\gamma\cdots$ shared by $P_2,P_3$ in the first of Figure \ref{3a2bCase3}. The remainder $\cdots$ of the vertex has angle sum $2\pi-\beta-\gamma=\frac{7}{12}\pi$. The only combination of angles with sum $\frac{7}{12}\pi$ is single $\beta$. This contradicts the fact that any single angle outside $P_2,P_3$ must be an $a^2$-angle. Similar contradiction happens to $f=18$ and $22$.

For $f=20$, we also get the values for the five angles, and find that the vertex $\beta\gamma\cdots$ can be $\beta^2\gamma$ or $\beta\gamma\epsilon^4$. The vertex $\beta^2\gamma$ has the same contradiction as the case of $f=16$. The vertex $\beta\gamma\epsilon^4$ implies that there is a vertex of degree $6$. By \eqref{vcountf}, we find $v_4=v_6=1$ and $v_5=v_7=v_8=\cdots=0$. This is combinatorially impossible by \cite[Theorems 6]{yan}.

\subsection{Cases 3.2 and 3.3}

The neighborhood tiling for Case 3.2 is the second of Figure \ref{3a2bCase3}. The neighborhood tiling has vertices $\alpha^3$, $\beta^2\delta$, $\beta\gamma\epsilon$, $\delta^2\epsilon$, and the angles are given by \eqref{eq_case32}. By $f\ge 16$, we find that $2\pi-2\epsilon$ is strictly smaller than all the angles. This implies that $\epsilon$ appears at most once at any vertex, and therefore Lemma \ref{klem} holds. We also note that $\gamma>0$ implies $f<60$. 

Next we find vertices $\alpha^a\beta^b\gamma^c\delta^d\epsilon$ with one $\epsilon$. The angle sum at the vertex gives
\begin{equation}\label{case2eq4}
\frac{2}{3}a
+\left(\frac{5}{6}-\frac{2}{f}\right)b
+\left(-\frac{1}{6}+\frac{10}{f}\right)c
+\left(\frac{1}{3}+\frac{4}{f}\right)d
=\frac{2}{3}+\frac{8}{f}.
\end{equation}
By $60>f\ge 16$, this implies $16a+17b+8d\le 28$. We substitute the finitely many values of $a,b,d$ satisfying the inequality into \eqref{case2eq4}, and find whether it is possible to have non-negative integer solution $c$. In addition to the existing vertices $\beta\gamma\epsilon,\delta^2\epsilon$, the other possibilities are $\alpha\gamma^c\epsilon$ ($c\ge 1$), $\gamma^c\delta\epsilon$ ($c\ge 1$), $\gamma^c\epsilon$ ($c\ge 2$). By Propositions \ref{beven} and \ref{klem}, we must have $c=2$ in these vertices. By $2\gamma+\epsilon=\left(1+\frac{12}{f}\right)\pi<2\pi$, $\gamma^2\epsilon$ is not a vertex.

If $\alpha\gamma^2\epsilon$ is a vertex, then the extra angle sum at the vertex implies $f=36$, which further implies
\[
\alpha=\frac{2}{3}\pi,\;
\beta=\frac{7}{9}\pi,\;
\gamma=\frac{1}{9}\pi,\;
\delta=\frac{4}{9}\pi,\;
\epsilon=\frac{10}{9}\pi.
\] 
Those vertex combinations satisfying Propositions \ref{beven} and \ref{klem} are
\[
\text{AVC}
=\{\alpha^3, \;
\beta\gamma\epsilon, \;
\beta^2\delta, \;
\delta^2\epsilon, \;
\alpha\beta\gamma\delta, \;
\alpha\gamma^2\epsilon\}.
\]
This implies that the vertex $\gamma^2\cdots$ shared by $P_2,P_3$ in the second of Figure \ref{3a2bCase3} is $\alpha\gamma^2\epsilon$. Then we find that $\alpha\epsilon$ is bounded by two $a$-edges, a contradiction. 

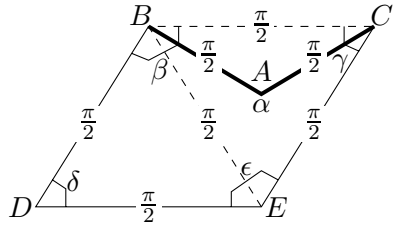
\begin{figure}[htp]
\centering
\begin{tikzpicture}[>=latex]

\draw[line width=1.5]
	(0,2.4) -- (1.5,1.5) -- (3,2.4);
		
\draw[dashed]
	(1.5,0) -- (0,2.4) -- (3,2.4);

\draw
	(0,2.4) -- (-1.5,0) -- (1.5,0) -- (3,2.4)
	(-1.1,0) -- (-1.1,0.24) -- (-1.28,0.36)
	(1.1,0) -- (1.1,0.2) -- (1.5,0.48) -- (1.72,0.36)
	(0.4,2.4) -- (0.4,2.16) -- (0,1.95) -- (-0.22,2.04)
	(2.6,2.4) -- (2.6,2.16) -- (2.8,2.07);

\node at (1.5,1.35) {\small $\alpha$};
\node at (0.15,1.8) {\small $\beta$};
\node at (2.55,1.9) {\small $\gamma$};
\node at (-1,0.35) {\small $\delta$};
\node at (1.3,0.5) {\small $\epsilon$};

\node[fill=white,inner sep=1] at (1.5,2.4) {\small $\frac{\pi}{2}$};
\node[fill=white,inner sep=1] at (-0.8,1.1) {\small $\frac{\pi}{2}$};
\node[fill=white,inner sep=1] at (0.8,1.1) {\small $\frac{\pi}{2}$};
\node[fill=white,inner sep=1] at (2.2,1.1) {\small $\frac{\pi}{2}$};
\node[fill=white,inner sep=1] at (0,0) {\small $\frac{\pi}{2}$};
\node[fill=white,inner sep=1] at (0.8,2) {\small $\frac{\pi}{2}$};
\node[fill=white,inner sep=1] at (2.2,2) {\small $\frac{\pi}{2}$};

\node at (1.5,1.8) {\small $A$};
\node at (-0.1,2.55) {\small $B$};
\node at (3.1,2.55) {\small $C$};
\node at (-1.7,0) {\small $D$};
\node at (1.7,0) {\small $E$};

\end{tikzpicture}
\caption{Impossible spherical pentagon.}
\label{case3fig1}
\end{figure}

If $\gamma^2\delta\epsilon$ is a vertex, then the extra angle sum at the vertex implies $f=24$, which further implies
\[
\alpha=\frac{2}{3}\pi,\;
\beta=\frac{3}{4}\pi,\;
\gamma=\frac{1}{4}\pi,\;
\delta=\frac{1}{2}\pi,\;
\epsilon=\pi.
\]
Substituting the angles into \eqref{eq10}, we get $L=2,M=0,N=0$, and then solve to get $\cos a=0$. Since $0<a<\pi$ (because $BD$ and $DE$ intersect only at $D$), we get $a=\frac{1}{2}\pi$. In Figure \ref{case3fig1}, by $BD=DE=a=\frac{1}{2}\pi$ and $\angle BDE=\frac{1}{2}\pi$, we know $\triangle BDE$ is an equilateral triangle with side length $\frac{1}{2}\pi$ and inner angle $\frac{1}{2}\pi$. Then by $BE=CE=\frac{1}{2}\pi$ and $\angle BEC=\angle CED-\angle BED=\pi-\frac{1}{2}\pi=\frac{1}{2}\pi$, we know $\triangle BCE$ is also an equilateral triangle with side length $\frac{1}{2}\pi$ and inner angle $\frac{1}{2}\pi$. Furthermore, by $\angle ABE=\angle ABD-\angle DBE=\frac{3}{4}\pi-\frac{1}{2}\pi=\frac{1}{4}\pi$ and $\angle ACE=\frac{1}{4}\pi$, the edges $AB$ and $AC$ evenly divide the angles $\angle CBE$ and $\angle BCE$. Therefore $A$ is the center of the equilateral triangle $\triangle BCE$, so that $\alpha=\frac{4}{3}\pi$, contradicting to $\alpha=\frac{2}{3}\pi$.

It remains to consider the case $\alpha\gamma^2\epsilon$ and $\gamma^2\delta\epsilon$ are not vertices. In this case, the only possible vertices in addition to the existing $\alpha^3$, $\beta\gamma\epsilon$, $\beta^2\delta$, $\delta^2\epsilon$ are $\alpha^a\beta^b\gamma^c\delta^d$. Similar to the argument in Section \ref{type2B}, to avoid contradicting Lemma \ref{klem}, each $\gamma$ needs to be combined with one $\beta$ into a chain $\beta\delta\cdots\delta\gamma$ bordered by two $b$-edges. This implies $b\ge c$. Now we know all vertices have $b\ge c$. In order to balance the equal total number $f$ of $\beta$ and $\gamma$, we must have $b=c$ at all the vertices. This contradicts the fact that $\beta^2\delta$ is a vertex in the second of Figure \ref{3a2bCase3}.

This completes the proof that Case 3.2 admits no tilings. 

Finally, we consider Case 3.2 given by the third of Figure \ref{3a2bCase3}. The neighborhood tiling has vertices $\alpha^3$, $\beta\gamma\delta$, $\gamma^2\epsilon$, $\delta\epsilon^2$. These are obtained from the vertices $\alpha^3$, $\beta\gamma\epsilon$, $\beta^2\delta$, $\delta^2\epsilon$ in the second of Figure \ref{3a2bCase3} by exchanging $\beta$ with $\gamma$ and exchanging $\delta$ with $\epsilon$ (the values in \eqref{eq_case33} are therefore the exchange of the values in \eqref{eq_case32}). In the argument for Case 3.2, the only other fact we used about the second of Figure \ref{3a2bCase3} is that $P_2,P_3$ sharing a vertex $\gamma^2\cdots$. The same exchange converts this to the vertex $\beta^2\cdots$ shared by $P_2,P_3$ in the third of Figure \ref{3a2bCase3}. Therefore after exchanging the angles, the proof for the impossibility of Case 3.2 becomes the proof for the impossibility of Case 3.3. We note that the proof also requires the exchange of angles in Lemma \ref{klem}.

\end{document}